\documentclass[11pt,a4paper]{amsart}

\usepackage[T1]{fontenc}
\usepackage[utf8]{inputenc}
\usepackage{lmodern}
\usepackage{amsmath,amssymb,amsthm,mathtools,mathrsfs}
\usepackage[margin=31mm]{geometry}
\usepackage{microtype}
\usepackage{xcolor}
\usepackage[colorlinks=true,citecolor=blue!55!black,
            linkcolor=blue!55!black,urlcolor=blue!55!black]{hyperref}

\newcommand{\Hh}{\mathbb H}
\newcommand{\Q}{\mathbb Q}
\newcommand{\R}{\mathbb R}
\newcommand{\Z}{\mathbb Z}
\newcommand{\e}{\mathrm e}
\newcommand{\ep}{\varepsilon}
\newcommand{\Li}{\operatorname{Li}_2}
\newcommand{\Cl}{\operatorname{Cl}_2}
\newcommand{\trans}{^{\mathsf{t}}}
\newcommand{\Norm}{\mathrm N}
\newcommand{\frakf}{\mathfrak f}
\newcommand{\Poch}[3]{(#1,#2;#3)_\infty}

\DeclareMathOperator{\im}{Im}
\DeclareMathOperator{\re}{Re}
\DeclareMathOperator{\Log}{Log}

\theoremstyle{plain}
\newtheorem{theorem}{Theorem}[section]
\newtheorem{proposition}[theorem]{Proposition}
\newtheorem{lemma}[theorem]{Lemma}
\newtheorem{corollary}[theorem]{Corollary}

\theoremstyle{definition}
\newtheorem{definition}[theorem]{Definition}
\newtheorem{example}[theorem]{Example}

\theoremstyle{remark}
\newtheorem{remark}[theorem]{Remark}

\numberwithin{equation}{section}

\title[Shintani invariants and cyclic quantum dilogarithms]
      {Shintani invariants and cyclic quantum dilogarithms}
\author{Bora Yalkinoglu}
\address{CNRS and IRMA, Universit\'e de Strasbourg}
\email{yalkinoglu@math.unistra.fr}
\date{}
\begin{document}

\begin{abstract}
We express general Shintani invariants of real quadratic number fields in terms of cyclic quantum dilogarithms, which appear in the asymptotic expansion of $q$-Pochhammer symbols. We also explain how a formula from 1995 of Sczech for Shintani invariants can naturally be understood in this framework.
\end{abstract}

\maketitle
\tableofcontents

\section{Introduction}

\noindent
In his seminal paper, Shintani attached a positive real number \(X(\mathfrak c) = X_1(\mathfrak c)X_2(\mathfrak c)\) to each
strict ray class \(\mathfrak c \in \mathrm{Cl}_K^+(\mathfrak f)\) in the strict ray class group of a real quadratic number field $K$, using the
double sine function \cite{shintani,yamamoto-kronecker}. Shintani conjectured that his invariants give rise to abelian extensions of \(K\) \cite{shintani-ray}, a prediction closely related to Stark's conjectures. In a breakthrough work, Radchenko and Wheeler proved algebraicity of Shintani invariants\cite{radchenko-wheeler}, the unit property and the action of the ray class group are given in \cite[Theorem~3]{stark-reciprocity}, solving Shintani's conjecture in general. In this work, which continues our announcement \cite{yalkinoglu-cqd-announcement}, we give a new expression for the Shintani invariants in terms of cyclic quantum dilogarithms.

\noindent
The starting point is an expression for the invariants in terms of the
\(q\)-Pochhammer symbol: Yamamoto found it in an example
\cite{yamamoto-factorization}, the author proved it in an imporatant special case \cite{yalkinoglu-note}, and Kopp gave the general
formula using his Shintani--Faddeev modular cocycle \cite{kopp}, which plays an important role in this paper as well. As \(q\) approaches a primitive $N$-th root of unity, the asymptotic expansion for the \(q\)-Pochhammer symbol contains the cyclic quantum dilogarithm, given by the
product \cite{bazhanov-reshetikhin}
\[
 D_{m/N}(x,y)=\prod_{j=1}^{N-1}\bigl(1-\e((j+x)\tfrac{m}{N}+y)\bigr)^{j/N},
 \qquad \e(z)=e^{2\pi iz}.
\]
Here \(m/N\) is in lowest terms and \(N>0\).


\noindent
For pedagogical reasons, we first consider the setting of our announcement
\cite{yalkinoglu-cqd-announcement}. Write \(\ep>1\) for the generator
of the totally positive units, \(1_{\frakf}\) for the identity ray
class of modulus \(\frakf\), and \(g\) for the order of \(\ep\bmod\frakf\).

\begin{theorem}\label{thm:intro-special}
Suppose \(\mathcal O_K=\Z[\ep]\), \(\frakf=(\alpha)\) is proper and
\(\alpha\gg0\). Write \(\alpha^{-1}\equiv x\ep+y\bmod\mathcal O_K\),
with \(0<x\le1\), \(0\le y<1\). Put \(a=\ep+\ep^{-1}\), and let
\(N_n\) be the reduced denominator of the Chebyshev fraction
\[
 t_n=\frac{T_{n-1}(a)}{T_n(a)},\qquad T_k(z+z^{-1})=z^k+z^{-k}.
\]
Let \(C_n(u,v)\) be the normalized quotient of
\(D_{t_n}(u,v)\) by \(D_{t_{n+g}}(u,v)\) in \eqref{def:Qn}. Then
\[
 \begin{aligned}
 C_n(y,x)&=X_1(1_{\frakf})\bigl(1+O(N_n^{-2})\bigr),\\
 C_n(x,y)&=X_2(1_{\frakf})\bigl(1+O(N_n^{-2})\bigr).
 \end{aligned}
\]
\end{theorem}

\noindent
Here a single sequence of fractions serves both invariants $X_1(\mathfrak c)$ and $X_2(\mathfrak c)$: only
\(x\) and \(y\) are interchanged. For an arbitrary ray class we keep
\((x,y)\) fixed and use a separate sequence for the other invariant. The main result of this paper is as follows.

\begin{theorem}\label{thm:intro-main}
Let \(\mathfrak c\) be any strict ray class of proper nonzero modulus
\(\frakf\). Choose \(\mathfrak b=\Z\omega+\Z\) in the narrow class of
\(\frakf\mathfrak a^{-1}\), where \(\mathfrak a\) represents
\(\mathfrak c\), with \(0<\omega'<1<\omega\).
Let \(x\omega+y\) represent a point of the Shintani orbit
\(\Phi_{\frakf,\mathfrak b}(\mathfrak c)\) of Definition~\ref{def:ray-orbit},
with \(0<x\le1,\ 0\le y<1\). Define
the integral matrix \(G\) by \(G(\omega,1)\trans=\ep^{-g}(\omega,1)\trans\).
For the explicit reduced fractions
\(m_n^+/N_n^+\to\omega\) and \(m_n^-/N_n^-\to\omega'\) of
\eqref{eq:v-paired}, the normalized cyclic quotients satisfy
\[
 \begin{aligned}
 C_G(m_n^+,N_n^+;x,y)&=X_1(\mathfrak c)\bigl(1+O((N_n^+)^{-2})\bigr),\\
 C_G(m_n^-,N_n^-;x,y)&=X_2(\mathfrak c)\bigl(1+O((N_n^-)^{-2})\bigr).
 \end{aligned}
\]
\end{theorem}

\noindent
As an application of our methods, we give a new approach to Sczech's approximation formula for Shintani invariants \cite{sczech} based on his function \(\Lambda(u,v;w)\). We express his rational values in terms of cyclic quantum dilogarithms. 

\begin{theorem}\label{thm:intro-sczech}
Let \(u_k,v_k\) be the traces and \(w_k\) the ratio of the two
basis elements at step \(k\) of the continued-fraction cycle for \(\mathfrak c\) in
Section~\ref{sec:sczech}, with \(0\le k<h\) and \(u_0\notin\Z\).
For the rational approximants of
Theorem~\ref{thm:sczech-global-corrected}, one has for all sufficiently large \(n\)
\[
 \sum_{k=0}^{h-1}\Lambda_*(u_k,v_k;w_{k,n})
 =\log C_H(m_n,N_n;\{-v_0\},u_0)
 =\log X_1(\mathfrak c)+O(N_n^{-2}).
\]
For the conjugate class the limit is \(\log X_2(\mathfrak c)\).
\end{theorem}
\subsection*{Outlook}
\noindent
There are several interesting directions for future research, starting with a better understanding of the connections between this paper, our recent work \cite{yalkinoglu-fermat} and the recent breakthrough of Radchenko and Wheeler \cite{radchenko-wheeler}. \\ \\
The cyclic quantum dilogarithm is naturally related to Fermat curves
\cite{faddeev-kashaev}, a comparison of the current techniques with the construction
of \cite{yalkinoglu-fermat} will appear in a follow-up work. \\ \\ 
Moreover, it's an interesting problem to relate the techniques of our paper to the case of complex cubic number fields, where the elliptic gamma function plays the role of the double sine function \cite{bergeron-charollois-garcia,morain}. It's also an interesting problem to understand links between the present work and recent \(p\)-adic breakthrough approaches to Stark's conjectures and explicit class field theory for real quadratic fields, such as \cite{dasgupta-kakde,darmon-vonk}.

\subsection*{Outline}
\noindent
Section~\ref{sec:shintani} attaches the Shintani orbit to a ray class and recalls earlier results from \cite{yalkinoglu-note,kopp}.
Section~\ref{sec:boundary} is the technical heart of our paper: it uses the asymptotics of the \(q\)-Pochhammer symbol at a root of unity to give an exact formula for the modulus of the cocycle in terms of the normalized cyclic quotient. The main result of that section, Theorem~\ref{thm:boundary}, is used to prove Theorem~\ref{thm:general-class} in Section~\ref{sec:general-classes}. In Section~\ref{sec:factor-relations} we give relations for Shintani invariants, inspired by the work of Yamamoto \cite{yamamoto-kronecker}. In Section~\ref{sec:sczech} we relate our results to Sczech's work \cite{sczech}.

\begingroup
\makeatletter
\let\@tocwrite\@gobbletwo
\makeatother
\section*{Notation and conventions}
\endgroup
\noindent
\emph{Arithmetic.}
\(K\) is real quadratic with ring of integers \(\mathcal O_K\).
We fix an embedding into \(\R\) and write \(a'=\sigma(a)\) for the other;
\(\operatorname{Tr}(a)=a+a'\), \(\Norm(a)=aa'\), and \(a\gg0\)
means \(a,a'>0\). For a nonzero integral ideal \(\mathfrak a\),
\(\Norm(\mathfrak a)=|\mathcal O_K/\mathfrak a|\).
Throughout \(0\ne\frakf\subsetneq\mathcal O_K\).
The strict ray class group \(\mathrm{Cl}_K^+(\frakf)\) is the group of
fractional ideals prime to \(\frakf\), modulo principal ideals \((a)\)
with \(a\gg0\) and \(a\equiv1\bmod\frakf\) at the primes dividing
\(\frakf\). Omitting the finite congruence gives the narrow class
group \(\mathrm{Cl}_K^+\). The group of totally positive units is
\(\langle\ep\rangle\), \(\ep>1\); those congruent to \(1\bmod\frakf\)
form \(\langle\ep^{g(\frakf)}\rangle\), where
\(g(\frakf)\) is the order of \(\ep\) in \((\mathcal O_K/\frakf)^\times\).

\noindent
\emph{Shintani invariants.}
For a strict ray class \(\mathfrak c\) modulo \(\frakf\), let
\[
 \zeta_{\frakf}(s,\mathfrak c)
 =\sum_{\substack{0\ne\mathfrak a\subseteq\mathcal O_K\\
                   \mathfrak a\in\mathfrak c}}
   \Norm(\mathfrak a)^{-s}
 \qquad (\Re(s)>1)
\]
be its partial zeta function. To specify the two classes entering
Shintani's formula, choose \(b\in\mathcal O_K\) with \(b\gg0\) and
\(b\equiv-1\pmod{\frakf}\), and put
\[
 \mathfrak J=[(b)]_{\frakf}\in\mathrm{Cl}_K^+(\frakf).
\]
This class depends only on \(\frakf\): the quotient of two such
choices is totally positive and congruent to \(1\) modulo \(\frakf\).
Moreover, \(\mathfrak J^2=1_{\frakf}\), since \(b^2\gg0\) and
\(b^2\equiv1\pmod{\frakf}\).

\noindent
The Shintani invariant \(X(\mathfrak c)>0\) is given by
\begin{equation}\label{eq:shintani-zeta}
 \log X(\mathfrak c)
 =\zeta'_{\frakf}(0,\mathfrak c\mathfrak J)
  -\zeta'_{\frakf}(0,\mathfrak c),
\end{equation}
where the prime denotes differentiation in \(s\) after analytic
continuation. Thus replacing \(\mathfrak c\) by
\(\mathfrak c\mathfrak J\) interchanges the two terms and replaces
\(X(\mathfrak c)\) by its reciprocal. We write
\(X(\mathfrak c)=X_1(\mathfrak c)X_2(\mathfrak c)\), with the
two factors given by Yamamoto's double-sine products at the chosen
embedding and its conjugate, respectively
\cite[\S5.2, (5.2.1)]{yamamoto-kronecker}. Their equivalent cocycle
formulas are given in \eqref{eq:kopp-endpoints}.

\noindent
\emph{Coordinates.}
A basis of an ideal is an ordered \(\Z\)-basis. Modulo
\(\Z\omega+\Z\) we write \(x\omega+y\) with \(0<x\le1,\ 0\le y<1\).
The two column conventions in the paper are
\(v=(x,y)\trans\) for these coordinates and \(r=(-y,x)\trans\)
in Kopp's cocycle. Write \(\{u\}=[u]=u-\lfloor u\rfloor\) and
\([u]_1=1-\{-u\}\). The brackets \([[b_0,b_1,\ldots]]\) denote a minus continued
fraction. For \(A=\left(\begin{smallmatrix}a&b\\c&d\end{smallmatrix}\right)\),
\(A\tau=(a\tau+b)/(c\tau+d)\), whereas \(Aw\)
denotes ordinary column multiplication.

\noindent
\emph{Limits and branches.}
There are two limits. In Section~\ref{sec:boundary} we fix the reduced
fraction \(m/n\), \(n>0\), and approach its root of unity along a radius:
\[
 \tau=\frac mn+\frac{it}{2\pi n^2}\in\Hh,\qquad
 q=\e(\tau)=\e(m/n)e^{-t/n^2},\qquad t\downarrow0\quad(t\in\R_{>0}).
\]
In Section~\ref{sec:general-classes} it is the fractions
\(m_j/n_j\) themselves that converge to a quadratic irrational, with
\(n_j\to\infty\). Our asymptotic notation is
\[
\begin{aligned}
 a_j=O(b_j)&:\quad |a_j|\le C |b_j|\text{ for all sufficiently large }j;\\
 a_j\asymp b_j&:\quad c b_j\le a_j\le C b_j\text{ for all sufficiently large }j\quad(a_j,b_j>0);\\
 a_j\sim b_j&:\quad a_j/b_j\longrightarrow1\quad(b_j\ne0\text{ eventually}).
\end{aligned}
\]
Here \(c,C>0\) are independent of \(j\), but may depend on the fixed
class, basis and matrix. For \(t\downarrow0\) the same inequalities
hold for all sufficiently small \(t>0\). We use principal branches of
\(\Log\) and \(\Li\), taking powers in a cyclic product one factor at
a time; \(\log\) on \(\R_{>0}\) is the real logarithm.

\section{Shintani invariants and \texorpdfstring{\(q\)}{q}-Pochhammer functions}
\label{sec:shintani}

\noindent
Write \([\mathfrak a]_{\frakf}\) for the strict ray class of an ideal prime
to \(\frakf\), and \([\mathfrak a]^+\) for its narrow ideal class. The map
\[
 \pi_{\frakf}:\mathrm{Cl}_K^+(\frakf)\longrightarrow\mathrm{Cl}_K^+,
 \qquad [\mathfrak a]_{\frakf}\longmapsto[\mathfrak a]^+
\]
forgets the congruence modulo \(\frakf\). Fix a narrow class \(\mathfrak h\).
We describe the ray classes above \(\mathfrak h\) as finite orbits on a
single lattice quotient. Choose a fractional ideal \(\mathfrak b\) with
\([\mathfrak b]^+=[\frakf]^+\mathfrak h^{-1}\). The reason for this
choice is that for every \(\mathfrak a\) in \(\mathfrak h\), the lattice
\(L=\frakf\mathfrak a^{-1}\) can be scaled to \(\mathfrak b\) by a
totally positive element. Our starting point on \(K/L\) is simply \(1+L\).

\begin{definition}\label{def:ray-orbit}
A \emph{Shintani orbit} of modulus \(\frakf\) on \(K/\mathfrak b\) is a set
\[
 \langle\ep\rangle z=\{\ep^jz:j\in\Z\},\qquad
 z\in K/\mathfrak b,\qquad
 \operatorname{Ann}(z):=\{a\in\mathcal O_K:az=0\}=\frakf .
\]
Write \(\operatorname{Sh}_{\frakf}(\mathfrak b)\) for the set of these orbits.
For an integral representative \(\mathfrak a\) prime to \(\frakf\), choose
\(\eta\gg0\) such that \(\eta\frakf\mathfrak a^{-1}=\mathfrak b\).
The Shintani orbit of its ray class is
\begin{equation}\label{eq:ray-torsion-map}
 \begin{aligned}
 \Phi_{\frakf,\mathfrak b}:\pi_{\frakf}^{-1}(\mathfrak h)
   &\longrightarrow\operatorname{Sh}_{\frakf}(\mathfrak b),\\
 [\mathfrak a]_{\frakf}&\longmapsto\langle\ep\rangle(\eta+\mathfrak b).
 \end{aligned}
\end{equation}
\end{definition}

Multiplication by \(\eta\) sends \(1+L\) to \(\eta+\mathfrak b\).
We take its orbit because any other positive scaling differs by a power
of \(\ep\). This is the fixed-lattice form of the ray class parametrization
in \cite[\S2.1]{yamamoto-kronecker} and \cite[Theorem~3.14]{kopp}.

\begin{lemma}\label{lem:ray-orbit}
The map \(\Phi_{\frakf,\mathfrak b}\) is a bijection. Every Shintani
orbit has \(g(\frakf)\) points. Changing \(\mathfrak b\) to
\(\lambda\mathfrak b\), \(\lambda\gg0\), carries each orbit to the
orbit of the same ray class by the isomorphism
\[
 K/\mathfrak b\longrightarrow K/(\lambda\mathfrak b),\qquad
 z+\mathfrak b\longmapsto\lambda z+\lambda\mathfrak b .
\]
\end{lemma}

\begin{proof}
Since \(\mathfrak a\) is prime to \(\frakf\),
\(\operatorname{Ann}(1+L)=L\cap\mathcal O_K=\frakf\).
Thus \(\ep^j(1+L)=1+L\) exactly when \(\ep^j\equiv1\bmod\frakf\),
which gives the period. If \(\delta\mathfrak a\) is another integral
representative of the ray class, take \(\delta\gg0\) with
\(\delta\equiv1\bmod\frakf\) at the primes of \(\frakf\).
Then \(\delta-1\in\frakf\mathfrak a^{-1}\): at primes outside
\(\frakf\) this follows from the integrality of both ideals.
The new scaling \(\delta\eta\) therefore gives the same point modulo
\(\mathfrak b\).

Conversely, given \(z\) with \(\operatorname{Ann}(z)=\frakf\), lift it to
\(\eta\gg0\) by adding a sufficiently large positive integer in
\(\mathfrak b\), and set \(\mathfrak a=(\eta)\frakf\mathfrak b^{-1}\).
This ideal is integral since \(\frakf z=0\). At each
\(\mathfrak p\mid\frakf\), the equality
\(\eta^{-1}\mathfrak b\cap\mathcal O_K=\frakf\) gives
\(v_{\mathfrak p}(\eta^{-1}\mathfrak b)=v_{\mathfrak p}(\frakf)\);
hence \(\mathfrak a\) is prime to \(\frakf\).
Any other positive lift has ratio to \(\eta\) congruent to \(1\bmod\frakf\)
at every prime dividing \(\frakf\), since their difference lies in
\(\mathfrak b=\eta\frakf\mathfrak a^{-1}\).
Replacing \(\eta\) by \(\ep^j\eta\) does not change its ideal.
Thus the inverse sends the Shintani orbit of \(z\) to
\([\mathfrak a]_{\frakf}\). Replacing \(\eta\) by \(\lambda\eta\)
proves the assertion about rescaling.
\end{proof}

For the formulas of Section~\ref{sec:general-classes} we need coordinates
for the very same points. From now on, for a class \(\mathfrak c\) we take
\(\mathfrak h=\pi_{\frakf}(\mathfrak c)\) and use a \emph{reduced}
representative \(\mathfrak b=\Z\omega+\Z\), meaning \(0<\omega'<1<\omega\).
Such a representative exists by minus continued-fraction reduction
\cite[\S2.1]{yamamoto-kronecker}. The narrow class of \(\mathfrak b\)
is fixed; neither this lattice nor the first point of the Shintani orbit
is unique.

\begin{lemma}\label{lem:orbit-coordinates}
Let \(\mathfrak a\) be an integral ideal prime to \(\frakf\) representing
\(\mathfrak c\), and choose a totally positive
generator \(\eta\) of \(\mathfrak a\mathfrak b\frakf^{-1}\).
One point of \(\Phi_{\frakf,\mathfrak b}(\mathfrak c)\) is represented
by \(x\omega+y\), where
\begin{equation}\label{eq:ray-torsion-coordinates}
 x=\left[\frac{\eta-\eta'}{\omega-\omega'}\right]_1,\qquad
 y=\left\{\frac{\omega\eta'-\omega'\eta}{\omega-\omega'}\right\}.
\end{equation}
Here \(a'\) denotes conjugation, and \([\,\cdot\,]_1,\{\,\cdot\,\}\)
take the representatives modulo \(\Z\) in \((0,1]\) and \([0,1)\).
Let \(H\) be the integral matrix defined by
\(H(\omega,1)\trans=\ep(\omega,1)\trans\). Starting at \(v_0=(x,y)\trans\),
the successive coordinate columns are
\[
 v_{j+1}\equiv H\trans v_j\pmod{\Z^2},\qquad 0\le j<g(\frakf).
\]
The columns \(v_0,\ldots,v_{g(\frakf)-1}\) represent the Shintani orbit;
\(v_{g(\frakf)}=v_0\) when the same coordinate representatives are used.
\end{lemma}

\begin{proof}
Write \(\eta=s\omega+t\) with \(s,t\in\Q\). Solving this equation and
its conjugate gives \eqref{eq:ray-torsion-coordinates} after reduction
modulo \(\Z\). Since \(\ep\mathfrak b=\mathfrak b\) and
\(\Norm(\ep)=1\), the matrix \(H\) belongs to \(\mathrm{SL}_2(\Z)\).
Multiplication of \(s\omega+t\) by \(\ep\) acts on the row \((s,t)\)
by right multiplication by \(H\), hence on its transpose by \(H\trans\).
Lemma~\ref{lem:ray-orbit} gives the first return after \(g(\frakf)\) steps.
\end{proof}

For principal ideals \(\frakf=(\alpha)\) and \(\mathfrak a=(\gamma)\)
with \(\alpha,\gamma\gg0\) and \(\mathfrak a+\frakf=\mathcal O_K\), take
\(\mathfrak b=\mathcal O_K\). Then \(\eta=\gamma/\alpha\). In particular,
for \(\frakf=(N)\), \(N\in\Z_{>1}\), the identity class starts at
\(1/N\), with coordinates \((1,1/N)\) in a basis \((\omega,1)\).

\begin{example}
Take \(K=\Q(\sqrt5)\), \(\ep=(3+\sqrt5)/2\), and
\(\omega=\ep\). Then \(H\trans(x,y)\trans=(3x+y,-x)\trans\).
For the identity class the two examples of
\cite[Examples~2.1-2.2]{yalkinoglu-note} give the following
Shintani orbits (each row returns to its first entry):
\[
\begin{array}{c@{\qquad}l}
 \frakf & \multicolumn{1}{l}{(x_j,y_j)\text{ in cyclic order}}\\[2pt]
 (4) & \tfrac14(4,1),\ \tfrac14(1,0),\ \tfrac14(3,3)\\[2pt]
 (4-\sqrt5) & \tfrac1{11}(2,1),\ \tfrac1{11}(7,9),\ \tfrac1{11}(8,4),\
               \tfrac1{11}(6,3),\ \tfrac1{11}(10,5).
\end{array}
\]
In the second row the initial point is
\(1/(4-\sqrt5)=(2\ep+1)/11\). Thus the two periods are \(3\) and \(5\).
\end{example}

\subsection{The product formula when \texorpdfstring{\(\mathcal O_K=\Z[\ep]\)}{the maximal order is Z[epsilon]}}
\label{subsec:unit-basis-input}

\noindent
Let \(\ep>1\) generate the totally positive units, write
\(\ep'=\ep^{-1}\) and \(a=\ep+\ep'\), and assume
\begin{equation}\label{eq:monogenic}
 \mathcal O_K=\Z+\Z\ep,\qquad
 U=\begin{pmatrix}a&-1\\1&0\end{pmatrix}.
\end{equation}
For \(\frakf=(\alpha)\), \(\alpha\gg0\), the Shintani orbit of the
identity ray class starts at the point
\begin{equation}\label{eq:initial-data}
 \alpha^{-1}\equiv x\ep+y\pmod{\mathcal O_K},\qquad
 (x,y)\in(0,1]\times[0,1).
\end{equation}
Since \(\frakf\) is proper, \((x,y)\notin\Z^2\). Put \(g=g(\frakf)\).
Multiplication by \(\ep\) gives
\((x_j,y_j)\trans\equiv (U\trans)^j(x,y)\trans\)
with the same representatives and least period \(g\).

\noindent
For \(\omega>0\) and \(0<z<1+\omega\), our double sine is
\[
 \mathcal S(\omega,z)=
 \exp\bigl(\zeta_2'(0,1+\omega-z\mid\omega,1)
          -\zeta_2'(0,z\mid\omega,1)\bigr).
\]
Here \(\zeta_2(s,z\mid\omega,1)\) continues
\(\sum_{m,n\ge0}(z+m\omega+n)^{-s}\) from \(\re s>2\), and the
prime denotes differentiation in \(s\). The two Shintani invariants are
\begin{equation}\label{eq:X12}
 X_i(1_{\frakf})=\prod_{j=1}^{g}\mathcal S(\theta_i,x_j\theta_i+y_j),
 \qquad (\theta_1,\theta_2)=(\ep,\ep').
\end{equation}
Write \((w;q)_\infty=\prod_{k\ge0}(1-wq^k)\) for \(|q|<1\).
For \(\tau\in\Hh\), abbreviate
\begin{equation}\label{eq:qpoch}
 \Poch{x}{y}{\tau}=(\e(x\tau+y);\e(\tau))_\infty=\varpi_{(-y,x)\trans}(\tau).
\end{equation}
The double-sine product formula, in the form of
\cite[Proposition~3.1]{yalkinoglu-note}, gives, for \(\theta>0\) and nonintegral
\((u,v)\in(0,1]\times[0,1)\),
\begin{equation}\label{eq:yamamoto-step}
 \mathcal S(\theta,u\theta+v)=
 \lim_{\substack{\tau\to\theta\\\tau\in\Hh}}
 \left|\frac{\Poch{u}{v}{\tau}}
 {\Poch{1-v}{u}{-1/\tau}}\right|.
\end{equation}
For the Chebyshev polynomials in the normalization
\(T_n(z+z^{-1})=z^n+z^{-n}\), put \(T_{-n}=T_n\) and
\begin{equation}\label{eq:tau-n}
 \tau_n=\frac{T_{n+1}(a)+i(\ep-\ep')}{T_n(a)}.
\end{equation}
These are the points of the geodesic used in \cite[Lemmas~3.2-3.3]{yalkinoglu-note}:
\(U^k\tau_n=\tau_{n+2k}\), \(\tau_n\to\ep\) and \(\tau_{-n}\to\ep'\)
as \(n\to\infty\), and \(\re\tau_{-n}=T_{n-1}(a)/T_n(a)\).

\begin{proposition}\label{prop:qpoch}
With the above data one has \cite[Theorem~3.2]{yalkinoglu-note}
\begin{align}
 X_1(1_{\frakf})&=\lim_{n\to\infty}
 \left|\frac{\Poch{x}{y}{\tau_{n-g}}}{\Poch{x}{y}{\tau_{n+g}}}\right|,
 \label{eq:q-X1}\\
 X_2(1_{\frakf})&=\lim_{n\to\infty}
 \left|\frac{\Poch{x}{y}{\tau_{-n-g}}}{\Poch{x}{y}{\tau_{-n+g}}}\right|.
 \label{eq:q-X2}
\end{align}
\end{proposition}

\noindent
The cancellation behind both formulas is the same:
\(1-y_j=x_{j-1}\) and \(x_j\equiv ax_{j-1}+y_{j-1}\pmod\Z\).
Thus the denominator of \eqref{eq:yamamoto-step} is
\(\Poch{x_{j-1}}{y_{j-1}}{U\tau}\). Multiplying over a period
at the arguments \(U^{g-j}\tau\) leaves just
\(\varpi_r(\tau)/\varpi_r(U^g\tau)\). We will use this same
cancellation for a reduced ideal cycle, with one matrix for each cone.

\subsection{Kopp's cocycle at the two real embeddings}

\noindent
For \(r=(-y,x)\trans\in\Q^2\setminus\Z^2\) and
\(\Gamma_r=\{A\in\mathrm{SL}_2(\Z):Ar-r\in\Z^2\}\), write
\begin{equation}\label{eq:sf-cocycle}
 \mathscr S_r(A;\tau)=\frac{\varpi_r(A\tau)}{\varpi_r(\tau)}
 \quad (A\in\Gamma_r,\ \tau\in\Hh),\qquad
 \varpi_r(\tau)=(\e(x\tau+y);\e(\tau))_\infty .
\end{equation}
This is Kopp's Shintani-Faddeev \emph{modular} cocycle (denoted by
the letter \emph{shin} in \cite[Definition~4.18]{kopp}); his Jacobi
cocycle \(\sigma_A(z,\tau)\), with independent elliptic variable \(z\),
is a different function. On \(\Hh\), cancellation gives
\begin{equation}\label{eq:analytic-cocycle}
 \mathscr S_r(AB;\tau)=\mathscr S_r(A;B\tau)\,\mathscr S_r(B;\tau).
\end{equation}

\begin{theorem}\label{thm:kopp-continuation}
For \(A=\left(\begin{smallmatrix}a&b\\c&d\end{smallmatrix}\right)\in\Gamma_r\),
the function \(\mathscr S_r(A;\cdot)\) extends meromorphically across
every real \(t\) with \(ct+d>0\)
\cite[Theorem~4.30]{kopp}.
\end{theorem}

\noindent
At \(m/n\) with \(cm/n+d>0\), the value
\(\mathscr S_r(A;m/n)\) is defined by continuation, not by substituting
\(q=\e(m/n)\) into the \(q\)-Pochhammer symbol. The theorem gives
neither a finite product nor, on its own, nonvanishing. Our
Theorem~\ref{thm:boundary} gives a positive finite formula for its modulus.

\noindent
The following is Kopp's reduced-cycle formula in Shintani's
two-embedding notation. We indicate precisely the passage from his
stated first-embedding formula to the second.

\begin{theorem}\label{thm:kopp-endpoints}
For \(\mathfrak c\in\mathrm{Cl}_K^+(\frakf)\), take
\((x,y)\) from Lemma~\ref{lem:orbit-coordinates} and put
\(r=(-y,x)\trans\), \(G=H^{-g(\frakf)}\). Then
\begin{equation}\label{eq:kopp-endpoints}
 X_1(\mathfrak c)=|\mathscr S_r(G;\omega)|,\qquad
 X_2(\mathfrak c)=|\mathscr S_r(G;\omega')|.
\end{equation}
Thus the auxiliary lattice determines the slopes and coordinates on
the right, but not the values on the left.
\end{theorem}

\begin{proof}
Kopp's formula \cite[Proposition~7.20]{kopp} gives
\(|\mathscr S_r(G^{-1};\omega)|=X_1(\mathfrak c)^{-1}\) for the
expanding eigenvalue. An integral shift of \(r\) does not affect
the value at a fixed point \cite[Proposition~4.35]{kopp}.
Inverting the matrix in \eqref{eq:analytic-cocycle} gives the first formula.

For the second, apply the same double-sine identity at the conjugate
slopes in Yamamoto's product \cite[Section~5.1]{yamamoto-kronecker}.
Explicitly let \(b_j\) be the digits of the reduced cycle,
\(V_j=\left(\begin{smallmatrix}b_j&-1\\1&0\end{smallmatrix}\right)\),
and \(L\) the period of the digits and the point together. For one
minimal period the product of the matrices is \(H\), since the
product of the successive reduced slopes is \(\ep\)
\cite[Section~2.1]{yamamoto-kronecker}. Each circuit thus advances
one point on the Shintani orbit. Lemma~\ref{lem:ray-orbit} gives the
first return after \(g(\frakf)\) circuits, so
\(V_0\cdots V_{L-1}=H^{g(\frakf)}=G^{-1}\). For
\(\tau_j=V_j\tau_{j+1}\) one has \(y_{j+1}=1-x_j\) and
\(x_{j+1}\equiv b_jx_j+y_j\bmod\Z\), whence
\[
 \Poch{1-y_{j+1}}{x_{j+1}}{-1/\tau_{j+1}}
     =\Poch{x_j}{y_j}{\tau_j}.
\]
In the product of \eqref{eq:yamamoto-step} over a period, all
intermediate terms cancel. The end quotient is
\(\varpi_r(\tau_L)/\varpi_r(G^{-1}\tau_L)=\mathscr S_r(G;G^{-1}\tau_L)\).
Let \(\tau_L\to\omega'\) in \(\Hh\). Each \(\tau_j\) tends to the
corresponding positive conjugate slope, so the product of the
double sines is \(X_2(\mathfrak c)\); \(G_{21}\omega'+G_{22}=\ep^{g(\frakf)}>0\) gives
the continued modulus on the right. This proves \eqref{eq:kopp-endpoints}.
\end{proof}

\noindent
For \(\mathfrak c=1_{\frakf}\) in Subsection~\ref{subsec:unit-basis-input},
\(G=U^{-g}\) and \((\omega,\omega')=(\ep,\ep')\):
\begin{equation}\label{eq:special-kopp-endpoints}
 X_1(1_{\frakf})=|\mathscr S_r(U^{-g};\ep)|,\qquad
 X_2(1_{\frakf})=|\mathscr S_r(U^{-g};\ep')|.
\end{equation}
These are also the two limits in Proposition~\ref{prop:qpoch}, since
\(U^{-g}\tau_{\pm n+g}=\tau_{\pm n-g}\).

\section{Cyclic quantum dilogarithm}
\label{sec:boundary}

\noindent
At a rational point the two \(q\)-Pochhammer symbols in
\(\mathscr S_r(A;\tau)\) may both tend to zero. Kopp's theorem gives
meromorphic continuation of their quotient; it does not give its value.
Theorem~\ref{thm:boundary} computes the modulus at every reduced
\(m/n\) with \(A_{21}m+A_{22}n>0\), including the cases in which a
cyclic quantum dilogarithm vanishes.

\subsection{Finite products and their normalization}

\noindent
A \emph{cusp} in this section means a reduced fraction \(m/n\),
represented by the primitive column \(w=(m,n)\trans\) with
\(n>0\). Let \(\tau\) tend vertically to the cusp \(m/n\) by setting
\[
 \tau(t)=\frac{m}{n}+\frac{it}{2\pi n^2}, \qquad t\downarrow0.
\]
Then \(q=\e(\tau(t))=\e(m/n)e^{-t/n^2}\) tends radially to the
root of unity \(\e(m/n)\).  For real \(x,y\), put
\begin{equation}\label{eq:D-def}
 D_{m/n}(x,y)
 =\prod_{j=1}^{n-1}
   \left(1-\e\left(\frac{(j+x)m}{n}+y\right)\right)^{j/n}.
\end{equation}
Complex powers in \eqref{eq:D-def} use the principal logarithm when a
complex value is required; a vanishing factor makes the product zero.
All main formulas use the modulus, which is branch independent.  We set
\(D_{m/1}(x,y)=1\).

\noindent
From now on fix \(r=(-y,x)\trans\in\Q^2\setminus\Z^2\),
\(0\le x\le1\). Let
\[
 A=\begin{pmatrix}\alpha&\beta\\c&d\end{pmatrix}\in\Gamma_r,
 \qquad Aw=(M,p)\trans,
\]
so that
\[
 M=\alpha m+\beta n,\qquad p=cm+dn.
\]
The pair \(M,p\) is again coprime.  We call the cusp \emph{admissible
for \(A\)} if \(p>0\), equivalently \(c(m/n)+d>0\).
Thus both denominators \(n,p\) are positive, and \(M/p\) is the
reduced representative of the transformed cusp.  We use only such
cusps.  Since \(A\in\Gamma_r\),
\begin{equation}\label{eq:rho-invariance}
 xM+yp\equiv xm+yn\pmod{\Z}.
\end{equation}
The common phase \(\rho=\e(xm+yn)\) is \emph{resonant} if
\(\rho=1\), and \emph{nonresonant} otherwise.  Since \(\gcd(m,n)=1\), the numbers
\[
 \e\!\left((j+x)\frac{m}{n}+y\right),\qquad 0\le j<n,
\]
are distinct, and each has \(n\)-th power \(\e(xm+yn)=\rho\).
They are therefore precisely the \(n\)-th roots of \(\rho\).
In particular, exactly one of them equals \(1\) if \(\rho=1\),
and none does otherwise.  Thus a factor with phase \(1\) occurs in the full
residue decomposition precisely at resonance. The next lemma makes this precise.

\begin{lemma}\label{lem:zero-match}
Suppose \(\rho=1\). Let \(0\le k<n\) and \(0\le\ell<p\) be the
unique integers for which
\begin{equation}\label{eq:zero-indices}
 (k+x)\frac{m}{n}+y\in\Z,\qquad
 (\ell+x)\frac{M}{p}+y\in\Z.
\end{equation}
Then
\[
 0<\frac{k+x}{n}=\frac{\ell+x}{p}<1.
\]
Moreover, \(D_{m/n}(x,y)=0\) if and only if \(k\ne0\), and
\(D_{M/p}(x,y)=0\) if and only if \(\ell\ne0\).
\end{lemma}

\begin{proof}
The preceding observation gives a unique index at each fraction,
since both have phase \(\rho=1\). Put \(s=(k+x)/n\). Then
\[
 s\binom{m}{n}\equiv\binom{-y}{x}=r\pmod{\Z^2}.
\]
We have \(0\le s\le1\), and either endpoint would give
\(r\in\Z^2\), which is excluded. Applying \(A\) and using
\(Ar\equiv r\pmod{\Z^2}\), we obtain
\[
 s\binom{M}{p}\equiv\binom{-y}{x}\pmod{\Z^2}.
\]
Thus \(ps-x\) is an integer in \([0,p)\) and \(sM+y\in\Z\).
By uniqueness, \(ps-x=\ell\), proving the equality.
The vanishing criteria follow because each cyclic product omits
the index \(0\).
\end{proof}
The following normalization gives a finite positive quotient even at
indices where the quotient in our announcement
\cite[Theorem~3.1]{yalkinoglu-cqd-announcement} is undefined, and is one of our new key definitions.
\begin{definition}
\label{def:corrected-quotient}
Fix \(A\in\Gamma_r\), where \(r=(-y,x)\trans\).
Let \(w=(m,n)\trans\) be a primitive integer column, and write
\(Aw=(M,p)\trans\), with \(n,p>0\). Put \(\rho=\e(xm+yn)\) and define
\begin{equation}\label{eq:intro-endpoint}
 E_r(w)=
 \begin{cases}
 |1-\rho|^{x/n}|D_{m/n}(x,y)|,&\rho\ne1,\\
 n^{-k/n}D_{m/n}^{\langle k\rangle}(x,y),&\rho=1 ,
 \end{cases}
 \qquad C_A(m,n;x,y)=\dfrac{E_r(w)}{E_r(Aw)} .
\end{equation}
Here \(D^{\langle k\rangle}\) is the product of the absolute values in
\eqref{eq:D-def} with \(j=k\) omitted (no omission if \(k=0\)).
Thus the \emph{normalized CQD quotient} is
\begin{align}
 C_A(m,n;x,y)
 &=|1-\rho|^{\,x(1/n-1/p)}
   \left|\frac{D_{m/n}(x,y)}{D_{M/p}(x,y)}\right|
      &&(\rho\ne1),\label{eq:Q-nonres}\\
 C_A(m,n;x,y)
 &=n^{-k/n}p^{\ell/p}
   \frac{D_{m/n}^{\langle k\rangle}(x,y)}
        {D_{M/p}^{\langle\ell\rangle}(x,y)}
      &&(\rho=1).\label{eq:Q-res}
\end{align}

\end{definition}

\noindent
These formulas include \(x=0,1\). All factors are positive algebraic,
since \(|1-\xi|^2=(1-\xi)(1-\xi^{-1})\) for a root of unity \(\xi\).

\begin{remark}\label{rem:matched-zero}
Suppose \(\rho=1\). To see where the factor
\(n^{-k/n}p^{\ell/p}\) in \eqref{eq:Q-res} comes from, put
\[
 P_{a/b}(t)=\prod_{j=1}^{b-1}
 \left|1-e^{-t/b}\e\!\left((j+x)\frac ab+y\right)\right|^{j/b},
 \qquad t>0.
\]
At \(m/n\), the term with index \(k\) contributes
\((1-e^{-t/n})^{k/n}\sim(t/n)^{k/n}\) as \(t\downarrow0\).
If \(k=0\), this contribution is \(1\), since the product starts
at \(j=1\). All other factors tend to nonzero values.
The same calculation at \(M/p\) gives \((t/p)^{\ell/p}\).

By Lemma~\ref{lem:zero-match},
\(k/n-\ell/p=x(1/p-1/n)\). Hence
\begin{equation}\label{eq:matched-zero-limit}
 C_A(m,n;x,y)=\lim_{t\downarrow0}
 t^{x(1/n-1/p)}\,\frac{P_{m/n}(t)}{P_{M/p}(t)}.
\end{equation}
The power of \(t\) cancels, leaving the constant
\(n^{-k/n}p^{\ell/p}\).

The two cyclic products need not even vanish together. For
\(A=\left(\begin{smallmatrix}3&-2\\-1&1\end{smallmatrix}\right)\),
\((x,y)=(1/2,0)\) and \(m/n=2/3\), one has \(M/p=0/1\):
the numerator vanishes at \(j=1\), whereas the denominator is
the empty product.
\end{remark}

\begin{lemma}\label{lem:finite-cocycle}\label{prop:finite-cocycle}
Let \(A,B\in\Gamma_r\), and write
\[
 w=\binom{m}{n},\qquad Bw=\binom{m'}{n'},\qquad
 ABw=\binom{m''}{n''},
\]
where \(w\) is primitive and \(n,n',n''>0\). Then
\begin{equation}\label{eq:finite-cocycle}
 C_{AB}(m,n;x,y)=C_A(m',n';x,y)\,C_B(m,n;x,y).
\end{equation}
\end{lemma}

\begin{proof}
The factor \(E_r(Bw)\) cancels in the product on the right.
\end{proof}

\subsection{Root-of-unity asymptotics}

\noindent
The asymptotic behaviour of the $q$-Pochhammer symbol at roots of unity is due to Bazhanov and Reshetikhin
\cite[Proposition~3.2, (3.10)-(3.11)]{bazhanov-reshetikhin};
see also \cite[Proposition~4.33]{kopp}. We need it also when the
argument tends to the unit circle.

\begin{proposition}\label{lem:BR}
Let \(\zeta\) be a primitive \(n\)-th root of unity and put
\(q=\zeta e^{-t/n^2}\). We use the following branch:
\[
 D_\zeta(w)=\prod_{j=1}^{n-1}(1-\zeta^jw)^j,\qquad
 D_\zeta(w)^{-1/n}
 :=\exp\!\left(-\sum_{j=1}^{n-1}\frac jn\Log(1-\zeta^jw)\right).
\]
As \(t\to0\) with \(|\arg t|\le\pi/2-\delta\) for fixed \(\delta>0\),
\begin{equation}\label{eq:BR-asymptotic}
 (w;q)_\infty
 =(1-w^n)^{1/2}\exp\!\left(-\frac{\Li(w^n)}t\right)
 D_\zeta(w)^{-1/n}\bigl(1+O(|t|)\bigr).
\end{equation}
This holds locally uniformly for \(|w|<1\). It also holds along
\(w=w_0e^{-\gamma t/n}\) if \(\gamma\ge0\) is fixed, \(|w_0|=1\),
and \(w_0^n\ne1\). The logarithms and square root in these formulas
are continued from \(w=0\), where \(\Log 1=0\) and \(\sqrt1=1\).
\end{proposition}

\begin{proof}
For fixed \(|w|<1\) the factor \(d(w)\) in
\cite[(3.10)]{bazhanov-reshetikhin} is
\((1-w^n)^{(n-1)/(2n)}D_\zeta(w)^{-1/n}\).  Substitution in
(3.11) there gives \eqref{eq:BR-asymptotic}.

For the moving argument, put \(T=t/n\), \(u_j=w_0\zeta^j\),
\(s_j=\gamma+j/n\). Grouping the \(q\)-Pochhammer factors by residues
modulo \(n\) gives factors \((u_j e^{-s_jT};e^{-T})_\infty\).
None of the \(u_j\) is \(1\).  For such \(u\) and fixed \(s\ge0\),
Euler-Maclaurin summation along the ray of \(T\) gives
\begin{equation}\label{eq:noncritical-block}
 \log(u e^{-sT};e^{-T})_\infty
 =-\frac{\Li(u)}{T}
  +\left(\frac12-s\right)\Log(1-u)+O(|T|).
\end{equation}
Indeed \(f_u(z)=\Log(1-u e^{-z})\) is regular at the origin and
\(\int_0^{\infty e^{i\arg T}} f_u(z)\,dz=-\Li(u)\).
On a closed subsector of \(\re T>0\), \(f_u'\) and \(f_u''\)
are bounded near \(0\) and decay exponentially at infinity,
uniformly along all its rays. The shifted endpoint contributes
\((1/2-s)f_u(0)\); the second-derivative remainder is \(O(|T|)\).
The same estimates are uniform for \(u\) in compact subsets of
the open disk and hence give the asserted local and sectorial
uniformity there as well.

Summing \eqref{eq:noncritical-block}, use
\(\sum_j\Li(u_j)=\Li(w_0^n)/n\) and
\(\sum_j\Log(1-u_j)=\Log(1-w_0^n)\), with \(0\le j<n\)
and the branches continued from the disk. The sum is
\[
 -\frac{\Li(w_0^n)}t+
 \left(\frac12-\gamma\right)\Log(1-w_0^n)
 -\sum_{j=1}^{n-1}\frac jn\Log(1-\zeta^jw_0)+O(|t|).
\]
Taylor expansion of \(\Li(w(t)^n)\) at \(w_0^n\ne1\) gives precisely
the same logarithm on the right of \eqref{eq:BR-asymptotic}.
\end{proof}

\noindent
For \(\zeta=\e(m/n)\), \(w=\e(xm/n+y)\) and \(w^n\ne1\),
\(D_\zeta(w)^{1/n}=D_{m/n}(x,y)\) with these branches.
The imaginary part of the dilogarithm on the unit circle is the
Clausen function
\[
 \Cl(\phi):=\sum_{j=1}^\infty\frac{\sin(j\phi)}{j^2}
 =\im\Li(e^{i\phi}),
\]
which plays a central role in our following key result.
\begin{theorem}
\label{thm:boundary}
Let \(r=(-y,x)\trans\in\Q^2\setminus\Z^2\), \(0\le x\le1\),
and \(A=\left(\begin{smallmatrix}\alpha&\beta\\c&d\end{smallmatrix}\right)
\in\Gamma_r\). At every admissible cusp \(m/n\) the meromorphic
continuation of \(\mathscr S_r(A;\tau)\) is regular and nonzero, with modulus
\begin{equation}\label{eq:boundary-algebraic}
 \left|\mathscr S_r\left(A;\frac mn\right)\right|
 = C_A(m,n;x,y)\,
   \exp\!\left(\frac{c}{2\pi np}\Cl\bigl(2\pi(xm+yn)\bigr)\right),
 \qquad p=cm+dn .
\end{equation}
At resonance the exponential is \(1\).  This is a formula for the
modulus; it does not prescribe a phase of the continued value.
\end{theorem}

\begin{proof}
Write \(\lambda=c/(2\pi i np)\).  The image of the radial path
satisfies the exact identity
\begin{equation}\label{eq:t-tilde}
 A\tau(t)=\frac Mp-\frac{\widetilde t}{2\pi i p^2},
 \qquad \widetilde t=\frac{t}{1-\lambda t}.
\end{equation}
Thus the second cusp uses \(\widetilde t\), not \(t\).
It lies in a fixed closed subsector of the right half-plane for
small \(t>0\).  We evaluate the two Pochhammer moduli there
using the same local calculation at each cusp.

If \(\rho\ne1\), the \(n\)-th power of the elliptic argument
\(\e(x\tau(t)+y)\) equals \(\rho e^{-xt/n}\).  Since
\(\Li(\rho e^{-a t})=\Li(\rho)+at\Log(1-\rho)+O(t^2)\),
Proposition~\ref{lem:BR} gives
\begin{equation}\label{eq:boundary-formula}
 |\varpi_r(\tau(t))|
 =\left|\exp\!\left(-\frac{\Li(\rho)}t\right)\right|
   \frac{|1-\rho|^{1/2}}{E_r(w)}(1+o(1)).
\end{equation}
The same formula holds at \(Aw\) with \(\widetilde t\), because
the phase is invariant by \eqref{eq:rho-invariance}.  In the
quotient, the factors \(|1-\rho|^{1/2}\) cancel.  The ratio of the
two exponential moduli is
\begin{equation}\label{eq:anomaly}
 \exp\!\left(\re\!\left[\left(\frac1t-\frac1{\widetilde t}\right)
               \Li(\rho)\right]\right)
 =\exp\!\left(\frac{c}{2\pi np}\Cl\bigl(2\pi(xm+yn)\bigr)\right).
\end{equation}
This proves the asserted radial modulus away from resonance.

At resonance there is a critical factor to which Proposition~\ref{lem:BR}
does not apply.  Put \(T=t/n\),
\(a_j=\e((j+x)m/n+y)\) and \(s_j=(j+x)/n\). The full product
splits as \(\prod_{j=0}^{n-1}(a_j e^{-s_jT};e^{-T})_\infty\).
For \(j\ne k\), use \eqref{eq:noncritical-block}; for \(j=k\)
we have \(a_k=1\), \(s_k=s\), and
\begin{equation}\label{eq:qgamma-limit}
 (e^{-sT};e^{-T})_\infty
 \sim\frac{\sqrt{2\pi}}{\Gamma(s)}
      T^{1/2-s}\exp\!\left(-\frac{\pi^2}{6T}\right).
\end{equation}
On the positive axis this is the usual \(q\)-gamma asymptotic
\cite[Section~1.10]{gasper-rahman}.  To see the sector uniformity
needed at the transformed cusp, write its logarithm, for \(h>1\), as
\[
 -\frac1{2\pi i}\int_{(h)}
   \Gamma(z)\zeta(z+1)\zeta(z,s)T^{-z}\,dz .
\]
Moving to \(\re z=-\eta\), \(0<\eta<1\), the residues at \(1,0\)
give \(-\pi^2/(6T)+(1/2-s)\Log T+\log(\sqrt{2\pi}/\Gamma(s))\).
Stirling's estimate bounds the remainder by \(O(|T|^\eta)\) on
closed subsectors of \(\re T>0\), for the fixed \(s>0\) in question.

Now \(\sum_{j=0}^{n-1}\Li(a_j)=\Li(1)/n\), while
\(\prod_{j\ne k}|1-a_j|=n\), since the \(a_j\) run through the
\(n\)-th roots of unity. The noncritical constant is therefore
\[
 \prod_{j\ne k}|1-a_j|^{1/2-(j+x)/n}
 =\frac{n^{1/2-x/n}}{D_{m/n}^{\langle k\rangle}(x,y)}.
\]
Multiplying by the power \(|t/n|^{1/2-s}\) from the critical
block leaves \(|t|^{1/2-s}/E_r(w)\), since \(s-x/n=k/n\).
Altogether the single-cusp modulus is
\[
 |\varpi_r(\tau(t))|
 =\left|\exp\!\left(-\frac{\Li(1)}t\right)\right|
   \frac{\sqrt{2\pi}}{\Gamma(s)}
   \frac{|t|^{1/2-s}}{E_r(w)}(1+o(1)).
\]
Lemma~\ref{lem:zero-match} gives the very same \(s\) at \(Aw\):
the gamma factors cancel and \(|\widetilde t/t|^{1/2-s}\to1\).
The ratio of the leading exponential moduli is \(1\), since
\(\lambda\) is purely imaginary and \(\Li(1)\) real.  What remains
is \(E_r(w)/E_r(Aw)\), as claimed.

By Theorem~\ref{thm:kopp-continuation}, \(c(m/n)+d>0\) puts this
rational point in the domain of meromorphic continuation. The radial
modulus we have calculated is finite and strictly positive, so the
Laurent order is zero and the limit is the modulus of the continued value.
\end{proof}

\begin{remark}\label{rem:kopp-correction}
For \(\im z>0\), Kopp's Jacobi cocycle with \(z\) fixed is
\[
 \sigma_A(z,\tau)=
 \frac{(\e(z/(c\tau+d));\e(A\tau))_\infty}{(\e(z);\e(\tau))_\infty}.
\]
Write \((m',p)\trans=A(m,n)\trans\), \(n,p>0\), and \(W=\e(nz)\).
The exponential below is missing from \cite[Proposition~4.34,
(4.17), arXiv v3]{kopp}:
\begin{equation}\label{eq:complex-anomaly}
 \sigma_A(z,m/n)=
 \frac{D_{\e(m'/p)}(\e(nz/p))^{-1/p}}
      {D_{\e(m/n)}(\e(z))^{-1/n}}\,
 \exp\!\left(\frac{c\,\Li(W)}{2\pi i np}
                 +\frac{cz}{p}\Log(1-W)\right).
\end{equation}
Indeed, for \(\widetilde z=z/(c\tau(t)+d)\),
\(\e(p\widetilde z)=W\exp((cz/p)\widetilde t)\) and
\(1/\widetilde t=1/t-c/(2\pi i np)\) by \eqref{eq:t-tilde}.
Thus the square roots in \eqref{eq:BR-asymptotic} cancel in the limit,
but Taylor expansion of the two dilogarithms gives the displayed
exponential, rather than \(1\) as in \cite[(4.18)]{kopp}.
The evaluation here has \(z\) fixed; it is not the modular cocycle
evaluation of Theorem~\ref{thm:boundary}.
\end{remark}

\noindent
Theorem~\ref{thm:boundary} uses \(Ar\equiv r\bmod\Z^2\).
When we change the ideal basis, the numerator instead uses a vector
congruent to \(Br\). We shall need the following formula for this quotient.

\begin{lemma}\label{lem:cqd-transfer}
Let \(B=\left(\begin{smallmatrix}a&b\\c&d\end{smallmatrix}\right)
\in\mathrm{SL}_2(\Z)\), \(s=(-y',x')\trans\equiv Br\bmod\Z^2\),
\(0\le x'\le1\). If \(w=(m,n)\trans\), \(Bw=(m',p)\trans\) and
\(n,p>0\), then the quotient below continues to \(m/n\), with modulus
\begin{equation}\label{eq:cqd-transfer-boundary}
 \left|\frac{\varpi_s(B\tau)}{\varpi_r(\tau)}\right|_{\tau=m/n}
 =\frac{E_r(w)}{E_s(Bw)}
  \exp\!\left(\frac{c}{2\pi np}\Cl(2\pi(xm+yn))\right).
\end{equation}
It is holomorphic and nonzero at each irrational real \(\theta\)
with \(c\theta+d>0\).
\end{lemma}

\begin{proof}
The congruence \(s\equiv Br\) gives \(x'm'+y'p\equiv xm+yn\bmod\Z\).
If this common value is nonintegral, \eqref{eq:boundary-formula} at
the two cusps proves \eqref{eq:cqd-transfer-boundary}. Otherwise, put
\(h=(m(k+x)+ny)/n\in\Z\), where \(k\) is given by
\eqref{eq:zero-indices}. Then the zero at the second cusp has index
\[
 \ell=p(k+x)/n-x'=ch+dk+(dx-cy-x')\in\Z:
\]
indeed \(0<(k+x)/n<1\) gives \(0\le\ell<p\), and
\(m'(k+x)/n+y'=ah+bk+(bx-ay+y')\in\Z\).
Thus \((k+x)/n=(\ell+x')/p\); the gamma factors cancel exactly as in
Theorem~\ref{thm:boundary}.

The quotient is the Jacobi cocycle with integral shift \(s-Br\).
By \cite[Remark~4.20 and Theorems~4.23, 4.29]{kopp} it continues
meromorphically across \(c\theta+d>0\). A zero or pole there requires
\(x\theta+y\in\Z\theta+\Z\), impossible for irrational \(\theta\)
and \(r\in\Q^2\setminus\Z^2\). At a rational point the strictly
positive finite radial modulus rules out a zero or pole.
\end{proof}

\subsection{Some properties of the normalized cyclic quotient}

\noindent
The following proposition explains the important distribution relation for the normalized cyclic quantum dilogarithm quotient \(C_A(m,n;x,y)\).

\begin{proposition}\label{rem:distribution}
Let \(q\ge1\) be an integer and let
\(r=(-y,x)\trans\in\Q^2\setminus\Z^2\), \(0\le x\le1\).
For \(0\le a,b<q\) choose the division representatives
\[
 x_b=(x+b)/q,\qquad y_a=(y+a)/q,\qquad
 s_{a,b}=(-y_a,x_b)\trans.
\]
Suppose \(A\in\mathrm{SL}_2(\Z)\) satisfies
\[
 A\equiv I\pmod q,\qquad (A-I)r\in q\Z^2 .
\]
Equivalently, \(A\) fixes every \(s_{a,b}\) modulo \(\Z^2\).
At every reduced \(m/n\) with \(n>0\) and
\(p=A_{21}m+A_{22}n>0\), one has
\begin{equation}\label{eq:finite-distribution}
 C_A(m,n;x,y)
 =\prod_{a,b=0}^{q-1}C_A(m,n;x_b,y_a),
\end{equation}
including all resonant cases.
\end{proposition}

\begin{proof}
The division points are the classes \((r+v)/q\) for
\(v\in\Z^2/q\Z^2\). Fixing the point with \(v=0\) and then fixing
all their differences gives the two stated congruences.

\noindent
At \(w=(m,n)\trans\) put \(z_j=\e((j+x)m/n+y)\).
Write \(d_n(z)=|1-z|\) for \(z\ne1\) and \(d_n(1)=1/n\).
Then in both cases of Definition~\ref{def:corrected-quotient},
\[
 E_r(w)=\prod_{j=0}^{n-1}d_n(z_j)^{(j+x)/n}.
\]
Indeed, the \(z_j\) are the \(n\)-th roots of \(\rho\): if \(\rho\ne1\)
the product of the \(|1-z_j|\) is \(|1-\rho|\); if \(z_k=1\), the
product with \(j=k\) omitted is \(n\).

\noindent
Multiply this expression over the division points. For fixed \(j,b\),
varying \(a\) gives all the \(q\)-th roots of \(z_l\), where
\(l\equiv qj+b\pmod n\), \(0\le l<n\).
Writing \(qj+b=l+in\) with \(0\le i<q\), their total exponent for
this \(l\) is
\(\sum_{i=0}^{q-1}(l+in+x)/(qn)=(l+x)/n+(q-1)/2\).
For \(z\ne1\), \(\prod_{\xi^q=z}d_n(\xi)=d_n(z)\); for \(z=1\)
it equals \(q/n=q\,d_n(1)\). Consequently, with \(d=(q-1)/2\),
\[
 \frac{\prod_{a,b}E_{s_{a,b}}(w)}{E_r(w)}
 =\begin{cases}
 |1-\rho|^d,&\rho\ne1,\\
 q^{d+(k+x)/n},&z_k=1.
 \end{cases}
\]
The right side is unchanged at \(Aw\): the phase \(\rho\) is the same,
and at resonance Lemma~\ref{lem:zero-match} gives
\((k+x)/n=(\ell+x)/p\). Dividing the two identities gives
\eqref{eq:finite-distribution}.
\end{proof}

\noindent
The corresponding distribution identity for the \(q\)-Pochhammer
symbol is the scalar case \(qI_2\) of
\cite[Proposition~4.45]{kopp}: this matrix leaves \(\tau\) unchanged
and sends \(s\) to \(qs\). For the analogous double-sine
multiplication formula see \cite[Theorem~2.1(b)]{kurokawa-koyama}.

\noindent
Reflection \(z\mapsto1+\omega-z\) inverts the double sine.
For the normalized cyclic quantum dilogarithm quotients the corresponding operation is
\((x,y)\mapsto(1-x,-y)\).

\begin{proposition}\label{prop:cqd-inversion}
For \(A\in\Gamma_r\), at every admissible cusp, including resonance,
\begin{equation}\label{eq:finite-inversion}
 C_A(m,n;x,y)C_A(m,n;1-x,-y)=1 .
\end{equation}
The two endpoint representatives give
\(C_A(m,n;1,y)=C_A(m,n;0,y)\) whenever \(y\notin\Z\) and
\(A\in\Gamma_{(-y,0)\trans}\).
Consequently \(C_A(w;v)\), with \(v=(x,y)\) represented in
\([0,1)^2\), satisfies \(C_A(w;-v)=C_A(w;v)^{-1}\).
\end{proposition}

\begin{proof}
Reversing the residues by \(j\mapsto n-1-j\) gives, off resonance,
\[
 |D_{m/n}(x,y)D_{m/n}(1-x,-y)|=|1-\rho|^{1-1/n}.
\]
Hence the corresponding two \(E\)-factors have product
\(|1-\rho|\), the same at \(w\) and \(Aw\).
At resonance the reflected critical residue is \(n-1-k\).
The product of the nonzero cyclotomic factors is \(n\), so the
two omitted products have product \(n^{(n-1)/n}\); the powers
\(n^{-k/n}\) and \(n^{-(n-1-k)/n}\) cancel it exactly.
This proves the first identity, also when \(x=0,1\).

\noindent
Shifting the residues by one gives
\(E_{(-y,1)\trans}(w)/E_{(-y,0)\trans}(w)
=|1-\e(y)|\), both off resonance and with the zero omitted.
Here \(y\notin\Z\).  The same factor occurs at \(Aw\), proving
the second identity.  Periodicity in \(y\) completes the proof.
\end{proof}

\noindent
A change of coordinates sends a reduced fraction \(m/n\) to
\[
 \frac{am+bn}{cm+dn},
 \qquad B=\begin{pmatrix}a&b\\c&d\end{pmatrix}\in\mathrm{GL}_2(\Z).
\]
Thus fractions tending to \(\theta\) are sent to fractions tending
to \(B\theta=(a\theta+b)/(c\theta+d)\), provided \(c\theta+d\ne0\).
The next proposition shows that the normalized cyclic quotient has
the same limit in the two coordinate systems.

\noindent
For \(r=(-y,x)\trans\) and \(w=(m,n)\trans\), write
\(C_A^r(w)=C_A(m,n;x,y)\). We write \(w\to\theta\) when
\(\gcd(m,n)=1\), \(n>0\), and \(m/n\to\theta\).

\begin{proposition}\label{prop:cqd-covariance}
Let \(A\in\Gamma_r\) fix a real quadratic irrational \(\theta\),
with \(A_{21}\theta+A_{22}>0\). Let
\(B=\left(\begin{smallmatrix}a&b\\c&d\end{smallmatrix}\right)
\in\mathrm{GL}_2(\Z)\) with \(c\theta+d>0\). Put
\[
 A'=BAB^{-1},\qquad
 s=(-y',x')\trans\equiv Br\pmod{\Z^2},\qquad 0\le x',y'<1.
\]
Then
\begin{equation}\label{eq:cqd-transfer-coboundary-limit}
 \lim_{w\to\theta}\frac{C_{A'}^s(Bw)}{C_A^r(w)}=1.
\end{equation}
In particular, the limiting moduli satisfy
\begin{equation}\label{eq:cqd-covariance}
 |\mathscr S_s(A';B\theta)|=|\mathscr S_r(A;\theta)|.
\end{equation}
\end{proposition}

\begin{proof}
Suppose first that \(\det B=1\). By definition,
\(C_A^r(w)=E_r(w)/E_r(Aw)\). Since \(A'Bw=BAw\), we obtain
\begin{equation}\label{eq:cqd-transfer-coboundary}
 \frac{C_{A'}^s(Bw)}{C_A^r(w)}
 =
 \frac{E_r(Aw)/E_s(BAw)}{E_r(w)/E_s(Bw)}.
\end{equation}
The numerator and denominator on the right are the same expression
evaluated at \(Aw\) and at \(w\). Both slopes tend to \(\theta\):
indeed, \(A(m/n)\to A\theta=\theta\). By
Lemma~\ref{lem:cqd-transfer}, both expressions therefore tend to
\[
 \left|\frac{\varpi_s(B\tau)}{\varpi_r(\tau)}\right|_{\tau=\theta}
 \ne 0.
\]
Here the Clausen exponential in the lemma tends to \(1\), since its
exponent is \(O(n^{-2})\). This proves
\eqref{eq:cqd-transfer-coboundary-limit}. If
\(|m/n-\theta|=O(n^{-2})\), holomorphy at \(\theta\) gives the
sharper estimate \(1+O(n^{-2})\) for the ratio in
\eqref{eq:cqd-transfer-coboundary}.

\noindent
For \(\det B=-1\), write \(B=J(JB)\), where
\(J=\operatorname{diag}(-1,1)\) and \(\det(JB)=1\).
The reflection \(J\) sends \(m/n\) to \(-m/n\) and
\((x,y)\) to \((x,-y)\). It complex conjugates every cyclic
factor, leaving its modulus unchanged; if a zero is omitted,
its index is unchanged as well. Thus the reflection preserves
the finite quotient exactly, and the result follows from the
case \(\det(JB)=1\).

\noindent
Finally, Theorem~\ref{thm:boundary} identifies the two limits
with the cocycle moduli in \eqref{eq:cqd-covariance}.
\end{proof}

\noindent
For example, interchanging the two basis vectors acts on slopes by
\(\theta\mapsto1/\theta\). When \(\theta=\ep\) this is precisely
\(\ep'\): Theorem~\ref{thm:main-cyclic} uses this observation to compute
\emph{both} Shintani invariants with the same Chebyshev fractions.
Equality \eqref{eq:cqd-covariance} is the modulus of
\cite[Theorem~4.37]{kopp}; \eqref{eq:cqd-transfer-coboundary} explains
why the finite approximations give it.

\section{Shintani invariants and cyclic quantum dilogarithms}
\label{sec:general-classes}
\label{sec:main}

\noindent
If \(\mathcal O_K=\Z[\ep]\) and \(\frakf=(\alpha)\) with \(\alpha\gg0\),
one sequence of rational numbers gives both invariants of the
identity class: interchanging \(x\) and \(y\) interchanges the embeddings.
For an arbitrary ray class we fix \((x,y)\) and instead approximate
each embedding separately.

\subsection{The Chebyshev formula for the identity class}
\label{subsec:special-cqd}
\noindent
Suppose \(\mathcal O_K=\Z[\ep]\). As in the announcement
\cite{yalkinoglu-cqd-announcement}, we take consecutive Chebyshev values
at \(a=\ep+\ep^{-1}\): \(T_n(a)=\ep^n+\ep^{-n}\).

\begin{lemma}\label{lem:cheb-reduction}
Let \(h=\gcd(a,2)\). For \(n\ge1\) put
\begin{equation}\label{eq:mn-def}
 w_n=\binom{m_n}{N_n}=\frac1h\binom{T_{n-1}(a)}{T_n(a)},
 \qquad t_n=m_n/N_n.
\end{equation}
Then \(w_n\) is primitive, so \(t_n\) is in lowest terms, and
\begin{equation}\label{eq:U-vector}
 U^{-1}w_n=w_{n+1},\qquad
 0<t_n-\ep'=\frac{(\ep-\ep^{-1})\ep^{-n}}{T_n(a)}=O(N_n^{-2}).
\end{equation}
\end{lemma}

\begin{proof}
The recurrence \(T_{n+1}(a)=aT_n(a)-T_{n-1}(a)\), with
\(T_0(a)=2,T_1(a)=a\), preserves the gcd of successive terms and
gives the matrix identity. The estimate follows by subtracting
\(\ep^{-1}T_n(a)\) from \(T_{n-1}(a)\).
\end{proof}

Set \(A=U^{-g}\); then \(Aw_n=w_{n+g}\).
For \(r=(-y,x)\trans\), periodicity gives
\((U\trans)^g(x,y)\trans\equiv(x,y)\trans\).
The identity \(JB\trans=B^{-1}J\) for \(B\in\mathrm{SL}_2(\Z)\), with
\(J=\left(\begin{smallmatrix}0&-1\\1&0\end{smallmatrix}\right)\), shows
that \(A\in\Gamma_r\).
For any admissible rational pair \((u,v)\), write
\begin{equation}\label{def:Qn}
 C_n(u,v)=C_A(m_n,N_n;u,v).
\end{equation}
Here admissibility means \(0\le u\le1\),
\((-v,u)\trans\notin\Z^2\), and
\(A\in\Gamma_{(-v,u)\trans}\).

\begin{theorem}\label{thm:main-cyclic}
Let \(\frakf=(\alpha)\) with \(\alpha\gg0\). The Shintani orbit
\(\Phi_{\frakf,\mathcal O_K}(1_{\frakf})\) is represented by
\(1/\alpha\). With its coordinates \((x,y)\) from \eqref{eq:initial-data},
\begin{align}
 \log C_n(y,x)&=\log X_1(1_{\frakf})+O(N_n^{-2}),\label{eq:rate-X1}\\
 \log C_n(x,y)&=\log X_2(1_{\frakf})+O(N_n^{-2}).\label{eq:rate-X2}
\end{align}
The estimates include indices at which an ordinary cyclic product vanishes.
\end{theorem}

\begin{proof}
By Theorem~\ref{thm:kopp-endpoints}, the two targets are the
moduli of \(\mathscr S_r(A;\cdot)\) at \(\ep\) and \(\ep'\).
To use the same fractions for both invariants, take
\(R=\left(\begin{smallmatrix}0&1\\1&0\end{smallmatrix}\right)\);
then \(R\ep=1/\ep=\ep'\) and \(RAR^{-1}=A^{-1}\).
The change of basis in Proposition~\ref{prop:cqd-covariance} yields
the first equality below. The inverse cocycle at the fixed point,
followed by the negation identity of Proposition~\ref{prop:cqd-inversion}
in the limit, gives the second:
\[
 |\mathscr S_r(A;\ep)|
 =|\mathscr S_{Rr}(A^{-1};\ep')|
 =|\mathscr S_{-Rr}(A;\ep')|.
\]
Since \(-Rr=(-x,y)\trans\), the pair for the first embedding is
\((y,x)\). Its column is again fixed by \(A\) modulo \(\Z^2\);
both estimates now follow by the same argument at \(\ep'\).

\noindent
For either pair \((u,v)=(x,y)\) or \((y,x)\),
Theorem~\ref{thm:boundary} gives
\[
 \log C_n(u,v)
 =\log|\mathscr S_{(-v,u)\trans}(A;t_n)|
 -\frac{A_{21}\,\Cl(2\pi(um_n+vN_n))}{2\pi N_nN_{n+g}}.
\]
The first term has a Taylor expansion at \(\ep'\); use
\eqref{eq:U-vector} for its error. The second term is also
\(O(N_n^{-2})\), since the Clausen function is bounded and
\(N_{n+g}\asymp N_n\).
\end{proof}

\begin{remark}\label{rem:raw-rate}
The announcement \cite[Theorem~3.1]{yalkinoglu-cqd-announcement} used the ordinary
CQD quotients, without the normalization. On nonresonant subsequences
they indeed give
\begin{align}
 X_1(1_{\frakf})&=\lim_n
 \left|\frac{D_{t_n}(y,x)}{D_{t_{n+g}}(y,x)}\right|,
 \label{eq:raw-X1}\\
 X_2(1_{\frakf})&=\lim_n
 \left|\frac{D_{t_n}(x,y)}{D_{t_{n+g}}(x,y)}\right|.
 \label{eq:raw-X2}
\end{align}
Each limit is along any infinite subsequence nonresonant for its
indicated pair: \(um_n+vN_n\notin\Z\) for \((u,v)=(y,x)\) or
\((x,y)\), respectively. Indeed, at such indices
\[
 \log C_n(u,v)-
 \log\left|\frac{D_{t_n}(u,v)}{D_{t_{n+g}}(u,v)}\right|
 =u(1/N_n-1/N_{n+g})\log|1-\e(um_n+vN_n)|.
\]
The phases belong to a finite set, so this difference is \(O(N_n^{-1})\).
Since \(N_{n+g}/N_n\to\ep^g>1\), on a subsequence of fixed
phase \(\rho\ne1\) its magnitude is comparable to
\(N_n^{-1}\) if \(u\log|1-\rho|\ne0\); otherwise it vanishes.
\end{remark}

\subsection{The principal narrow class}
\label{subsec:principal-fibre}

\noindent
Suppose again that \(\mathcal O_K=\Z[\ep]\) and
\(\frakf=(\alpha)\), \(\alpha\gg0\). Every class in
\(\ker\pi_{\frakf}\) has an integral representative \((\gamma)\)
with \(\gamma\gg0\). We may therefore use \(\mathfrak b=\mathcal O_K\)
throughout this fibre. The fractions \(m_n/N_n\) are unchanged;
in place of \(1/\alpha\) we use \(\gamma/\alpha\).

\begin{theorem}\label{thm:principal-fibre}
Let \(\mathfrak c=[(\gamma)]_{\frakf}\in\ker\pi_{\frakf}\), with
\(\gamma\gg0\) integral and prime to \(\frakf\).
The Shintani orbit \(\Phi_{\frakf,\mathcal O_K}(\mathfrak c)\) is represented by
\(\gamma/\alpha\). Writing
\(\gamma/\alpha\equiv x\ep+y\pmod{\mathcal O_K}\) as in
Lemma~\ref{lem:orbit-coordinates}, we have
\begin{align}
 \log C_n(y,x)&=\log X_1(\mathfrak c)+O(N_n^{-2}),\label{eq:class-X1}\\
 \log C_n(x,y)&=\log X_2(\mathfrak c)+O(N_n^{-2}).\label{eq:class-X2}
\end{align}
\end{theorem}

\begin{proof}
In Lemma~\ref{lem:orbit-coordinates} we now have
\(\mathfrak a\mathfrak b\frakf^{-1}=(\gamma/\alpha)\).
The residues of its totally positive generators form the Shintani orbit
modulo \(\mathcal O_K\). The period is still \(g\),
so the proof of Theorem~\ref{thm:main-cyclic} applies at each point.
\end{proof}

\subsection{Finite approximations for an arbitrary ray class}
\label{subsec:general-data}

Let \(\mathfrak c\in\mathrm{Cl}_K^+(\frakf)\). Use the choice
of \(\mathfrak b=\Z\omega+\Z\) from Section~\ref{sec:shintani}.
The two eigenlines of multiplication by \(\ep\) give the two
sequences of fractions. To compute its matrix \(H\), write
\[
 \omega=[[\overline{b_0,\ldots,b_{\ell-1}}]],\qquad
 V_k=\begin{pmatrix}b_k&-1\\1&0\end{pmatrix},\qquad
 H=V_0\cdots V_{\ell-1},
\]
using the minimal minus continued-fraction period.
Thus \(H(\omega,1)\trans=\ep(\omega,1)\trans\) as in
Lemma~\ref{lem:orbit-coordinates}. Put \(g=g(\frakf)\) and \(G=H^{-g}\).
The matrix \(H\) has two uses here. Its transpose permutes the
coordinates of the finite Shintani orbit modulo \(\Z^2\).
To approximate the two slopes we instead apply opposite powers of \(H\)
to the integer column \(e_2=(0,1)\trans\), without reducing modulo \(\Z^2\):
\begin{equation}\label{eq:v-paired}
 v_n^\pm=\binom{m_n^\pm}{N_n^\pm}
 =\mp H^{\pm n}e_2.
\end{equation}

\begin{theorem}\label{thm:general-class}
For \(\mathfrak c\in\mathrm{Cl}_K^+(\frakf)\), let \((x,y)\) be the
coordinates of any point in the Shintani orbit \(\Phi_{\frakf,\mathfrak b}(\mathfrak c)\),
as in Lemma~\ref{lem:orbit-coordinates}. Then
\begin{align}
 \log C_G(m_n^+,N_n^+;x,y)
 &=\log X_1(\mathfrak c)+O((N_n^+)^{-2}),\label{eq:general-class-X1}\\
 \log C_G(m_n^-,N_n^-;x,y)
 &=\log X_2(\mathfrak c)+O((N_n^-)^{-2}).\label{eq:general-class-X2}
\end{align}
The estimates hold for all sufficiently large \(n\), whether or not
\(xm_n^\pm+yN_n^\pm\) is an integer.
\end{theorem}

\begin{proof}
First observe that Theorem~\ref{thm:boundary} applies to
\emph{any} sequence of reduced fractions \(t_j=m_j/n_j\to\theta\),
\(n_j>0\), if \(\theta\) is real quadratic, \(B\theta=\theta\),
\(B\in\Gamma_r\), \(r=(-y,x)\trans\in\Q^2\setminus\Z^2\), \(0\le x\le1\),
and \(B_{21}\theta+B_{22}>0\). It gives
\begin{equation}\label{eq:general-rm-rate}
 \log C_B(m_j,n_j;x,y)
 =\log|\mathscr S_r(B;\theta)|+O\bigl(|t_j-\theta|+n_j^{-2}\bigr).
\end{equation}
Indeed, the new denominator is
\(n_j(B_{21}t_j+B_{22})\asymp n_j\), so the bounded Clausen
term contributes \(O(n_j^{-2})\). Lemma~\ref{lem:cqd-transfer}
gives a regular nonzero cocycle at \(\theta\), hence a log modulus
varying by \(O(|t_j-\theta|)\).

\noindent
For our class, \(G\in\Gamma_r\) since
\(J(H\trans)^g=H^{-g}J\) and \((H\trans)^g(x,y)\trans\equiv(x,y)\trans\).
The two moduli in \eqref{eq:kopp-endpoints} are the desired invariants.
It remains to estimate the slopes.

\noindent
Put \(e_\theta=(\theta,1)\trans\). In the decomposition
\[
 \binom01=\frac{-\omega'e_\omega+\omega e_{\omega'}}{\omega-\omega'},
 \qquad He_\omega=\ep e_\omega,\quad He_{\omega'}=\ep^{-1}e_{\omega'},
\]
the coefficient of the growing vector in \(\mp H^{\pm n}(0,1)\trans\)
is positive in both cases. Hence \(N_n^\pm\asymp\ep^n\); after division
by \(N_n^\pm\), the decaying component contributes
\(O(\ep^{-n}/N_n^\pm)=O((N_n^\pm)^{-2})\) to the slope error.
Unimodularity gives primitivity, and \(Gv_n^\pm=v_{n\mp g}^\pm\)
gives eventual positivity also after transformation.
Applying \eqref{eq:general-rm-rate} at each endpoint proves the
two estimates.
\end{proof}

\begin{remark}\label{prop:primitive-seed}
The choice of the initial column \(e_2\) in \eqref{eq:v-paired}
is not essential. For any primitive \(v\in\Z^2\), take the sign of
each column \(H^{\pm n}v\) so that its denominator is positive for
all sufficiently large \(n\). Neither eigenline contains a nonzero
rational column, so the proof gives the same limits and the same
estimates in terms of the denominators. In particular, the Chebyshev fractions in
Theorem~\ref{thm:main-cyclic} arise by applying \(U^{-n}\) to
\(v=(a/h,2/h)\trans\).

\noindent
If the approximating fractions are free to choose, both Shintani
invariants can be computed without encountering a zero. Take
\[
 v=\begin{cases}
 (1,0)\trans,&x\notin\Z,\\
 (0,1)\trans,&x\in\Z.
 \end{cases}
\]
Then \((x,y)v\notin\Z\), since \((x,y)\notin\Z^2\). Use
\(G^{-j}v\) to approximate \(\omega\) and \(G^jv\) to approximate
\(\omega'\), choosing the sign of each column so that its
denominator is positive. By \eqref{eq:rho-invariance},
\[
 (x,y)G^k v\equiv(x,y)v\pmod{\Z}\qquad(k\in\Z).
\]
Thus every resulting column \((m,N)\trans\) satisfies
\(xm+yN\notin\Z\), and no cyclic factor vanishes. For a prescribed
sequence of fractions this need not hold; in that case
Definition~\ref{def:corrected-quotient} also covers the indices
at which a factor vanishes.
\end{remark}

\section{Relations for Shintani invariants}
\label{sec:factor-relations}

\noindent
When are \(X_1(\mathfrak c)\) and \(X_2(\mathfrak c)\) equal, and
when are they reciprocal? Yamamoto studied how these invariants
change when one changes the signature of a ray class: one invariant
is unchanged and the other is inverted
\cite[Theorem~5.2.3 and Remark~5.2.4(1)]{yamamoto-kronecker}.
We recover his identities from the symmetries of the cyclic
quotients. Conjugation, on the other hand, interchanges the two
real embeddings; in particular, \(X_1(\mathfrak c)=X_2(\mathfrak c)\)
whenever both the modulus and the class are fixed by conjugation.

\noindent
The cyclic formulas also give a direct way to compare the two
invariants. In Theorem~\ref{thm:principal-fibre} they are computed
from the same fractions, with \((x,y)\) and \((y,x)\) interchanged.
If \((y,x)\) belongs to the Shintani orbit of \((x,y)\), the two
invariants are equal; if \((-y,-x)\) belongs to that orbit, they
are reciprocal. These relations already hold for the geometric
means of the finite cyclic quotients over the orbit, before taking
limits. We conclude with explicit congruence criteria for ray
classes modulo \((\ell)\), where \(\ell\) is an odd rational prime
unramified in \(K\).

\subsection{Conjugation and the two signature classes}

\noindent
Let \(\mathfrak c=[\mathfrak a]_{\frakf}\). If we replace
\(\mathfrak a\) by \(\mu\mathfrak a\), where
\(\mu\equiv1\pmod{\frakf}\) and \(\mu\) is totally positive,
the strict ray class does not change. What happens to
\(X_1(\mathfrak c)\) and \(X_2(\mathfrak c)\) if \(\mu\) is
negative at just one of the two real embeddings?

\noindent
Choose \(\mu_1,\mu_2\in1+\frakf\) with
\(\mu_1<0<\mu_1'\) and \(\mu_2'<0<\mu_2\). Write
\[
 \mathfrak s_i=[(\mu_i)]_{\frakf},\qquad
 \mathfrak c\mathfrak s_i=[\mu_i\mathfrak a]_{\frakf}.
\]
Only the signs matter: two choices with the same signs have
totally positive quotient congruent to \(1\) modulo \(\frakf\).
In particular, \(\mathfrak s_i^2=1_{\frakf}\). Changing both
signs gives precisely the class used in \eqref{eq:shintani-zeta}:
\[
 \mathfrak s_1\mathfrak s_2=\mathfrak J,
\]
since \(-\mu_1\mu_2\) is totally positive and congruent to
\(-1\) modulo \(\frakf\).

\noindent
We also compare \(\mathfrak c\) with its conjugate
\(\mathfrak c^\sigma=[\mathfrak a^\sigma]_{\frakf^\sigma}\).

\begin{proposition}\label{prop:galois-exchange}
For every strict ray class \(\mathfrak c\),
\[
 X_1(\mathfrak c)
 =X_2(\mathfrak c^\sigma),\qquad
 X_2(\mathfrak c)
 =X_1(\mathfrak c^\sigma),
\]
and
\begin{align}
 X_1(\mathfrak c\mathfrak s_1)&=X_1(\mathfrak c),&
 X_2(\mathfrak c\mathfrak s_1)&=X_2(\mathfrak c)^{-1},
 \label{eq:signature-1}\\
 X_1(\mathfrak c\mathfrak s_2)&=X_1(\mathfrak c)^{-1},&
 X_2(\mathfrak c\mathfrak s_2)&=X_2(\mathfrak c).
 \label{eq:signature-2}
\end{align}
For a stable modulus \(\frakf^\sigma=\frakf\), it follows that
\[
\begin{aligned}
 \mathfrak c^\sigma=\mathfrak c
 &\ \Rightarrow\ X_1(\mathfrak c)=X_2(\mathfrak c),\\
 \mathfrak c^\sigma=\mathfrak c\mathfrak J
 &\ \Rightarrow\ X_1(\mathfrak c)=X_2(\mathfrak c)^{-1}.
\end{aligned}
\]
\end{proposition}

\begin{proof}
By Theorem~\ref{thm:kopp-endpoints} it suffices to compare two
cocycle values. We use two identities: replacing \(r\) by \(-r\)
inverts the cyclic quotient \eqref{eq:finite-inversion}; replacing the
ideal basis changes the column by the rule of
Proposition~\ref{prop:cqd-covariance}. First change one sign:
represent the original and translated classes by
\(\mathfrak a_0\) and \(\mathfrak a_1=\mu\mathfrak a_0\),
choosing \(\mu\equiv1\pmod{\frakf}\) integral with the prescribed signs.
For \(i=0,1\), take the reduced lattice
\(\mathfrak b_i=(\eta_i)\mathfrak a_i^{-1}\frakf
=\Z\omega_i+\Z\) with \(\eta_i\gg0\).  Then
\[
 \mathfrak b_1=\lambda\mathfrak b_0,\qquad
 \lambda=\frac{\eta_1}{\mu\eta_0},\qquad
 \binom{\omega_1}{1}=\lambda B\binom{\omega_0}{1},
 \quad B\in\mathrm{GL}_2(\Z).
\]
Since \(\omega_i>\omega_i'\) for both reduced bases,
\(\det B=\operatorname{sgn}\Norm(\lambda)=-1\), and
 \(B_{21}\omega_0+B_{22}=1/\lambda\). Conjugating gives the
corresponding equality at \(\omega_0'\). Put \(v_i=(x_i,y_i)\trans\) and
\(r_i=(-y_i,x_i)\trans\).
Since \(\eta_1/\lambda=\mu\eta_0\equiv\eta_0\pmod{\mathfrak b_0}\),
the columns and matrices satisfy
\[
 B\trans v_1\equiv v_0\pmod{\Z^2},\qquad
 r_1\equiv-Br_0\pmod{\Z^2},\qquad
 G_1=BG_0B^{-1}.
\]
At either embedding \(\theta=\omega_0,\omega_0'\) the new column is
\(\operatorname{sgn}(B_{21}\theta+B_{22})Br_0\)
by Proposition~\ref{prop:cqd-covariance}. The denominator has the
sign of \(\mu\) at that embedding. Thus the column is \(r_1\)
where \(\mu\) is negative, and \(-r_1\) where it is positive.
Proposition~\ref{prop:cqd-inversion} gives the reciprocal in the latter case.
This proves \eqref{eq:signature-1}-\eqref{eq:signature-2}.

\noindent
For conjugation, take \(\mathfrak a_1=\mathfrak a_0^\sigma\) and
\(\frakf_1=\frakf_0^\sigma\).  Now
\[
 \binom{\omega_1}{1}=\lambda B\binom{\omega_0'}{1},
 \qquad \lambda=\eta_1/\eta_0'\gg0.
\]
The reversed orientation again gives \(\det B=-1\) and
\(r_1\equiv-Br_0\), while conjugation replaces \(\ep\) by
\(\ep^{-1}\).  Thus \(G_1=BG_0^{-1}B^{-1}\), since
\(g(\frakf^\sigma)=g(\frakf)\).
Now \(B_{21}\omega_0'+B_{22}=1/\lambda>0\), also after conjugation.
Changing basis, inverting \(G_0\) at its fixed point, and negating
the column give
\[
 |\mathscr S_{r_1}(G_1;\omega_1)|
 =|\mathscr S_{-r_0}(G_0^{-1};\omega_0')|
 =|\mathscr S_{r_0}(G_0;\omega_0')|.
\]
This gives the conjugation identities. The last two assertions now
follow by taking \(\mathfrak c^\sigma=\mathfrak c\) or
\(\mathfrak c^\sigma=\mathfrak c\mathfrak J\).
\end{proof}

\noindent
In particular one can recover each \(X_i\) from the undivided
invariant \(X\) alone \cite[\S5.2]{yamamoto-kronecker}:
\begin{equation}\label{eq:yamamoto-factors-zeta}
 X_i(\mathfrak c)^2=X(\mathfrak c)\,X(\mathfrak c\mathfrak s_i).
\end{equation}
Indeed, substitute \(X=X_1X_2\) and use
\eqref{eq:signature-1}-\eqref{eq:signature-2}.

\begin{corollary}\label{cor:signature-consequences}
\label{cor:minus-one-triviality}
One has
\begin{align}
 X(\mathfrak c\mathfrak s_1)&=X_1(\mathfrak c)/X_2(\mathfrak c),&
 X(\mathfrak c\mathfrak s_2)&=X_2(\mathfrak c)/X_1(\mathfrak c),
 \label{eq:X-signature}\\
 X(\mathfrak c\mathfrak J)&=X(\mathfrak c)^{-1}.
 \label{eq:X-total-signature}
\end{align}
In particular,
\[
 X_1(\mathfrak c)=X_2(\mathfrak c)
 \iff X(\mathfrak c\mathfrak s_1)=1
 \iff X(\mathfrak c\mathfrak s_2)=1.
\]
\noindent
Fix a narrow ideal class \(\mathfrak h\). For each embedding, the
product over the ray classes above \(\mathfrak h\) is \(1\):
\[
 \prod_{\mathfrak c\in\pi_{\frakf}^{-1}(\mathfrak h)}
 X_i(\mathfrak c)=1 \qquad (i=1,2).
\]
The same holds for \(X=X_1X_2\), and hence also when the products
are taken over all strict ray classes modulo \(\frakf\).

\noindent
A unit \(u\equiv1\pmod{\frakf}\) of mixed signs gives a conclusion
for each ray class \(\mathfrak c\) separately:
\[
 u<0<u' \ \Rightarrow\ X_2(\mathfrak c)=1,\qquad
 u'<0<u \ \Rightarrow\ X_1(\mathfrak c)=1.
\]
If, for some \(j\),
\begin{equation}\label{eq:epsilon-minus-one}
 \ep^j\equiv-1\pmod{\frakf},
\end{equation}
then \(X_1(\mathfrak c)=X_2(\mathfrak c)=1\) for every
\(\mathfrak c\).
\end{corollary}

\begin{proof}
The identities for \(X\) follow from \(X=X_1X_2\) and
Proposition~\ref{prop:galois-exchange}.

\noindent
To see why the products are \(1\), pair \(\mathfrak c\) with
\(\mathfrak c\mathfrak J\). These two ray classes have the same
narrow class because \(\mathfrak J=[(b)]_{\frakf}\) with \(b\gg0\).
Since \(\mathfrak J=\mathfrak s_1\mathfrak s_2\), the two
signature identities give
\[
 X_i(\mathfrak c\mathfrak J)=X_i(\mathfrak c)^{-1}
 \qquad (i=1,2).
\]
Thus each pair contributes \(1\). If a class is paired with
itself, its positive value \(X_i(\mathfrak c)\) is already \(1\).

\noindent
Finally, if \(u<0<u'\), then
\(\mathfrak s_1=[(u)]_{\frakf}=1_{\frakf}\), since \(u\) is a unit.
The signature identity for \(\mathfrak s_1\) therefore gives
\(X_2(\mathfrak c)=X_2(\mathfrak c)^{-1}\). The case \(u'<0<u\)
is identical with the indices interchanged. Under
\eqref{eq:epsilon-minus-one}, the totally positive unit \(\ep^j\)
represents \(\mathfrak J\), so \(\mathfrak J=1_{\frakf}\) and
both invariants are \(1\).
\end{proof}

\noindent
\subsection{When the two cyclic products agree}
\label{subsec:arithmetic-symmetries}

In Theorem~\ref{thm:principal-fibre} the two limits use the same
fractions; only the coordinates differ. Put
\(v=(x,y)\trans,\ R=\left(\begin{smallmatrix}0&1\\1&0\end{smallmatrix}\right)\)
and \(P=U\trans\). With coordinates taken modulo \(\Z^2\), set
\[
 \mathcal O(v)=\{v,Pv,\ldots,P^{g-1}v\}\pmod{\Z^2}.
\]
This is the Shintani orbit \(\Phi_{\frakf,\mathcal O_K}(\mathfrak c)\)
in coordinates. We test whether it contains \(Rv\), \(-Rv\), or \(-v\).

\begin{proposition}\label{prop:factor-criteria}
\label{cor:orbit-criterion}
For a class \(\mathfrak c\) and column \(v\) as above,
\[
 \begin{aligned}
   Rv\in\mathcal O(v)&\quad\Rightarrow\quad
                   X_1(\mathfrak c)=X_2(\mathfrak c),\\
  -Rv\in\mathcal O(v)&\quad\Rightarrow\quad
                   X_1(\mathfrak c)X_2(\mathfrak c)=1,\\
   -v\in\mathcal O(v)&\quad\Rightarrow\quad
                   X_1(\mathfrak c)=X_2(\mathfrak c)=1 .
 \end{aligned}
\]
\end{proposition}

\begin{proof}
Write \(C_n(w)=C_n(u,v')\) for coordinates
\(w=(u,v')\trans\) represented in \([0,1)^2\);
\eqref{eq:finite-inversion} permits \(u=0\) in place of \(u=1\).
For \(w\not\equiv0\) with \(P^gw\equiv w\bmod\Z^2\), form
\[
 \widehat{C}_n(w)=\left(\prod_{j=0}^{g-1}C_n(P^jw)\right)^{1/g}.
\]
Moving \(w\) to \(Pw\) permutes the factors, while negating \(w\)
inverts them term by term by \eqref{eq:finite-inversion}.  Thus,
at every index, including resonance,
\[
 \widehat{C}_n(Pw)=\widehat{C}_n(w),\qquad
 \widehat{C}_n(-w)=1/\widehat{C}_n(w).
\]
If \(Rv\) or \(-Rv\) lies in \(\mathcal O(v)\), these identities
make the two means at \(v,Rv\) equal or reciprocal, respectively.
If \(-v\in\mathcal O(v)\), the positive mean at \(v\) is its own
reciprocal, hence \(1\); since \(RPR=P^{-1}\), the orbit of
\(Rv\) likewise contains \(-Rv\), so its mean is also \(1\).
These are exact identities for the finite means, not necessarily
for the individual quotients. To identify their limits, every
\(P^jv\) represents the same class,
and \(P^jRv=RP^{-j}v\) is a swap of one of its representatives.
Applying Theorem~\ref{thm:principal-fibre} to these finitely many
terms gives
\[
 \lim_n\widehat{C}_n(v)=X_2(\mathfrak c),\qquad
 \lim_n\widehat{C}_n(Rv)=X_1(\mathfrak c),
\]
with additive errors \(O(N_n^{-2})\) in their logarithms.
\end{proof}

For \(\omega=\ep\), the coordinate swap works because
\(R\ep=1/\ep=\ep'\) and \(\det R=-1\). The same argument applies
to an ideal \(\Z\omega+\Z\) whenever there is an integral matrix
with these two properties. The next criterion does not assume that
every ideal admits such a matrix.

\begin{proposition}\label{prop:reverser}
Let \(\omega,H,r=(-y,x)\trans\) be as in Theorem~\ref{thm:general-class}.
Suppose there are integers \(a,b,c,d\) such that
\[
 T=\begin{pmatrix}a&b\\c&d\end{pmatrix},\qquad
 ad-bc=-1,\qquad \frac{a\omega+b}{c\omega+d}=\omega',\qquad
 c\omega+d>0 .
\]
Compare \(-T\,r\) with the finite set \(\{H^j r\bmod\Z^2:j\in\Z\}\).
For any \(j\),
\[
\begin{aligned}
 -T\,r\equiv H^j r\pmod{\Z^2}
 &\quad\Rightarrow\quad X_1(\mathfrak c)=X_2(\mathfrak c),\\
 -T\,r\equiv-H^j r\pmod{\Z^2}
 &\quad\Rightarrow\quad X_1(\mathfrak c)=X_2(\mathfrak c)^{-1}.
\end{aligned}
\]
\end{proposition}

\begin{proof}
The point of \(T\) is that both invariants can now be evaluated
along the \emph{same} fractions \(m_n^-/N_n^-\to\omega'\). Writing
\(w_n^-=(m_n^-,N_n^-)\trans\), we claim that
\begin{equation}\label{eq:reverser-base}
 X_1(\mathfrak c)=\lim_n C_G^{-T\,r}(w_n^-),\qquad
 X_2(\mathfrak c)=\lim_n C_G^r(w_n^-).
\end{equation}
The second equality is Theorem~\ref{thm:general-class}. For the
first, conjugating the equation \(T\omega=\omega'\) gives
\(T\omega'=\omega\). Thus \(T\) exchanges the eigenlines of \(G\),
so \(TGT^{-1}=G^{-1}\). Changing basis and inverting at the fixed point gives
\[
 X_1(\mathfrak c)=|\mathscr S_{T\,r}(G^{-1};\omega')|
                 =|\mathscr S_{-T\,r}(G;\omega')|.
\]
Here we used Propositions~\ref{prop:cqd-covariance} and~\ref{prop:cqd-inversion};
\eqref{eq:general-rm-rate} gives the first limit.
Replacing \(r\) by \(H^jr\) leaves the limiting modulus unchanged:
\(H^j\) fixes \(\omega'\) with positive eigenvalue and commutes with
\(G\). Replacing \(r\) by \(-r\) inverts each finite quotient.
These are the two assertions. The roles of the embeddings can also
be reversed, since \((c\omega+d)(c\omega'+d)=-\det T=1\).
\end{proof}

\begin{lemma}\label{lem:negative-stabilizer}
Let \(B,D\in\Gamma_r\) fix the same real quadratic irrational
\(\theta\). If
\[
 B_{21}\theta+B_{22}>0
 \qquad\text{and}\qquad
 D_{21}\theta+D_{22}<0,
\]
then
\begin{equation}\label{cor:negative-stabilizer}
 |\mathscr S_r(B;\theta)|=1.
\end{equation}
In particular, in the notation of Theorem~\ref{thm:general-class},
if some \(D\in\Gamma_r\) fixes \(\omega\) and satisfies
\(D_{21}\omega+D_{22}<0\), then
\(X_1(\mathfrak c)=X_2(\mathfrak c)=1\).
\end{lemma}

\begin{proof}
Both \(B\) and \(D\) fix \(\theta\) and its conjugate, so they
are diagonal in the same eigenbasis and hence commute. Apply
Proposition~\ref{prop:cqd-covariance} with change of coordinates
\(-D\). This fixes \(\theta\), leaves \(B\) unchanged, and sends
\(r\) to \(-r\) modulo \(\Z^2\). Thus
\[
 |\mathscr S_r(B;\theta)|
 =|\mathscr S_{-r}(B;\theta)|
 =|\mathscr S_r(B;\theta)|^{-1},
\]
where the second equality follows from
Proposition~\ref{prop:cqd-inversion} by taking limits.
The modulus is positive, so it equals \(1\).

For the last assertion take \(B=G\). The eigenvalues of \(G\)
at \(\omega\) and \(\omega'\) are both positive. Those of \(D\)
have product \(\det D=1\), so if one is negative, both are.
Applying the first assertion at the two embeddings gives the
claim by \eqref{eq:kopp-endpoints}.
\end{proof}

\noindent
For an integer modulus, the equality of the two Shintani invariants
of the identity class and the one-variable cyclic dilogarithm formula
were announced in \cite[Theorem~3.2, equation~(3.3)]{yalkinoglu-cqd-announcement}.
The following corollary specifies the distinction between the two
limits: the normalized quotient is defined for every index, and
the ordinary quotient is taken over those \(n\) for which
\(u\nmid N_n\).

\begin{corollary}\label{cor:integer-conductor}
Suppose \(\mathcal O_K=\Z[\ep]\) and \(\frakf=(u)\) for an
integer \(u>1\). Write
\(D_{m/n}(z)=D_{m/n}(0,z)\). For the identity ray class,
\begin{equation}\label{eq:integer-one-variable}
 X_1(1_{\frakf})=X_2(1_{\frakf})
 =\lim_n C_n(0,1/u)
 =\lim_{\substack{n\to\infty\\u\nmid N_n}}
 \left|\frac{D_{t_n}(1/u)}{D_{t_{n+g}}(1/u)}\right|.
\end{equation}
\end{corollary}

\begin{proof}
The modulus and identity class are fixed by conjugation, so
Proposition~\ref{prop:galois-exchange} gives \(X_1(1_{\frakf})=X_2(1_{\frakf})\).
The normalized coordinates of \(1/u\) modulo \(\mathcal O_K\) are
\((1,1/u)\). Thus \(X_2(1_{\frakf})=\lim_n C_n(1,1/u)=\lim_n C_n(0,1/u)\)
by Theorem~\ref{thm:main-cyclic} and the endpoint identity of
Proposition~\ref{prop:cqd-inversion}. With first coordinate zero
the nonresonant correction is \(1\); resonance means \(u\mid N_n\),
where \(C_n\) uses the matched-zero prescription. Infinitely many
indices remain since \(\gcd(N_n,N_{n+1})=1\).
\end{proof}

\subsection{Equality and reciprocity modulo a rational prime}

\noindent
Let \(\ell\) be an odd rational prime unramified in \(K\).
A ray class modulo \((\ell)\) above the principal narrow class
can be written as \(\mathfrak c=[(\alpha)]_{(\ell)}\), with
\(\alpha\in\mathcal O_K\) totally positive and prime to \(\ell\).
Its conjugate is represented by \(\alpha'\). Thus we can compare
\(\mathfrak c\) with \(\mathfrak c^\sigma\) by comparing two
residues modulo \(\ell\).

\begin{proposition}\label{prop:prime-counts}\label{prop:prime-dichotomy}
Let \(\mathfrak c=[(\alpha)]_{(\ell)}\) be as above. Then
\[
\begin{aligned}
 \alpha'&\equiv \ep^j\alpha\pmod{\ell}
 &&\text{for some \(j\)}
 &\quad&\Longrightarrow\quad X_1(\mathfrak c)=X_2(\mathfrak c),\\
 \alpha'&\equiv-\ep^j\alpha\pmod{\ell}
 &&\text{for some \(j\)}
 &\quad&\Longrightarrow\quad X_1(\mathfrak c)X_2(\mathfrak c)=1.
\end{aligned}
\]
\end{proposition}

\begin{proof}
Two totally positive generators represent the same ray class
precisely when their residues differ by a power of \(\ep\).
Thus the first congruence says exactly that
\(\mathfrak c^\sigma=\mathfrak c\). Since the residue \(-1\)
represents \(\mathfrak J\), the second says exactly that
\(\mathfrak c^\sigma=\mathfrak c\mathfrak J\).
Both conclusions follow from Proposition~\ref{prop:galois-exchange}.
\end{proof}

\noindent
When \(\mathcal O_K=\Z[\ep]\), these are the two cyclic limits
of Theorem~\ref{thm:principal-fibre}:
\[
 X_1(\mathfrak c)=\lim_n C_n(y,x),\qquad
 X_2(\mathfrak c)=\lim_n C_n(x,y),
 \qquad
 \frac{\alpha}{\ell}\equiv x\ep+y\pmod{\mathcal O_K},
\]
where \(0<x\le1\) and \(0\le y<1\).

\noindent
How many classes satisfy these congruences? The following count
concerns the congruences themselves, not all possible equalities
between values of the invariants.

\begin{corollary}
Let \(g\) be the order of \(\ep\) modulo \(\ell\), and put
\[
 h=\begin{cases}
      \ell-1,&\ell\text{ splits in }K,\\
      \ell+1,&\ell\text{ is inert in }K.
    \end{cases}
\]
There are \((\ell-1)h/g\) ray classes above the principal narrow
class. Exactly \(\ell-1\) satisfy each congruence of
Proposition~\ref{prop:prime-counts}.

If \(g\) is odd, the two sets of classes are disjoint. Together
they exhaust the principal narrow class precisely when \(h=2g\).

If \(g\) is even, the two sets coincide. In this case the stronger
conclusion
\[
 X_1(\mathfrak c)=X_2(\mathfrak c)=1
\]
holds for every ray class modulo \((\ell)\), not only those
above the principal narrow class.
\end{corollary}

\begin{proof}
Put \(A=(\mathcal O_K/\ell\mathcal O_K)^\times\) and
\(E=\langle\ep\bmod\ell\rangle\). The ray classes above the
principal narrow class are parametrized by \(A/E\): every
invertible residue has a totally positive lift, and two residues
give the same class precisely when they differ by an element of \(E\).

The map \(a\mapsto a'/a\) on \(A\) has kernel
\(\mathbb F_\ell^\times\) and image the cyclic group of residues
of norm one. This group has order \(h\): in the split case
conjugation exchanges the two components of
\(\mathbb F_\ell^\times\times\mathbb F_\ell^\times\); in the
inert case \(a'/a=a^{\ell-1}\) in \(\mathbb F_{\ell^2}^\times\).
In particular, \(|A/E|=(\ell-1)h/g\).

On ray classes the map becomes
\[
 A/E\longrightarrow A/E,\qquad aE\longmapsto (a'/a)E.
\]
It is well defined since \(\ep'=\ep^{-1}\). Its image has order
\(h/g\), since \(E\) lies in the norm-one group. Each element of
the image therefore has exactly \(\ell-1\) preimages. The two
congruences ask for the preimages of \(E\) and \(-E\), respectively.

If \(g\) is odd, then \(-1\notin E\), so these are distinct
preimages and exhaust \(A/E\) exactly when \(h/g=2\).
If \(g\) is even, then \(\ep^{g/2}\equiv-1\pmod{\ell}\), since
\(-1\) is the unique element of order two in the norm-one group.
Thus \(E=-E\), and Corollary~\ref{cor:minus-one-triviality}
gives the assertion for every ray class. In particular, this
applies if \(g>\ell\), since then \(g=h=\ell+1\).
\end{proof}

\section{Sczech's formula for the Shintani invariants}
\label{sec:sczech}

\noindent
We want to explain how Sczech's formula \cite{sczech} for the Shintani invariants
using his function \(\Lambda(u,v;w)\) fits into our framework.
Sczech expresses \(\Lambda(u,v;w)\) at a rational number
\(w=a/c\) as two finite sums \cite[Theorem~1]{sczech}. We write each
in terms of a cyclic dilogarithm. For a ray class \(\mathfrak c\), Sczech sums one term for each
step of the continued fraction until the slope and the two traces
return to their initial values, the traces being taken modulo
\(\Z\). At rational approximants the intermediate terms cancel,
and the sum is the logarithm of a single normalized cyclic
quotient. Passing to the limit identifies Sczech's sum with
\(\log X_1(\mathfrak c)\); applying the same formula to the
conjugate class gives \(\log X_2(\mathfrak c)\).

\subsection{Sczech's function at rational arguments}

\noindent
For \(w>0\) and \(u,v\in\R\setminus\Z\), Sczech defines
\cite[(2.2), (3.1)]{sczech}
\[
 \begin{aligned}
 \frac{4\pi}{i}\Lambda(u,v;w)
 &=\sum_{j\ne0}\frac{w\,\e(jv)+w^{-1}\e(ju)}{j|j|}\\
 &\quad+\lim_{M,N\to\infty}
 \sum_{\substack{0<|m|\le M\\0<|n|\le N}}
 \frac{\operatorname{sgn}(wm+n)}{mn}\e(mu+nv).
 \end{aligned}
\]
The cutoffs tend independently to infinity, and \(\operatorname{sgn}(0)=0\). Put
\(\ell(z)=\log|2\sin(\pi z)|\) for \(z\notin\Z\), and
\(P_1(t)=\{t\}-\tfrac12\) for \(t\notin\Z\), \(P_1(t)=0\) otherwise.
The value at \(w=1\) and the distribution relation
\cite[(3.5), Lemma~1]{sczech} give, for positive integers \(a,c\),
\begin{equation}\label{eq:sczech-w1}
 \begin{aligned}
 \Lambda(u,v;1)&=P_1(u-v)\bigl(\ell(u)-\ell(v)\bigr),\\
 \Lambda(u,v;a/c)&=\sum_{j=0}^{a-1}\sum_{k=0}^{c-1}
 \Lambda\!\left(\frac{u+j}{a},\frac{v+k}{c};1\right).
 \end{aligned}
\end{equation}

\noindent
For coprime \(a,c>0\), put \(\Delta=cu-av\).
In our notation, Sczech's finite formula \cite[Theorem~1]{sczech} reads
\[
\begin{aligned}
 \Lambda(u,v;a/c)
 &=\sum_{j=0}^{a-1}
   P_1\!\left(\frac{cj+\Delta}{a}\right)
   \ell\!\left(\frac{u+j}{a}\right)\\
 &\quad+\sum_{k=0}^{c-1}
   P_1\!\left(\frac{ak-\Delta}{c}\right)
   \ell\!\left(\frac{v+k}{c}\right).
\end{aligned}
\]
The second sum is obtained from the first by interchanging
\((a,u)\) and \((c,v)\), which changes \(\Delta\) to \(-\Delta\).
It therefore suffices to express the first sum in terms of a cyclic
dilogarithm. We denote this expression by \(\mathcal F(a,c;u,\Delta)\). For coprime \(N,q>0\),
\(z\in\R\setminus\Z\) and \(\Delta\in\R\), write
\[
 h=\lfloor\Delta\rfloor,\qquad \delta=\Delta-h,\qquad
 m=q^{-1}\bmod N,\qquad \eta=(z-mh)/N,
\]
with \(0\le m<N\) and \(m=0\) if \(N=1\). Set
\begin{equation}\label{eq:F-local}
 \mathcal F(N,q;z,\Delta)=\log|D_{m/N}(0,\eta)|+
 \begin{cases}
   (\delta/N-\tfrac12)\ell(z),&0<\delta<1,\\
   \tfrac12\bigl(\ell(\eta)-\ell(z)\bigr),&\delta=0.
 \end{cases}
\end{equation}

\begin{theorem}\label{thm:sczech-cqd}
For \(a,c>0\) coprime and \(u,v\in\R\setminus\Z\), put \(\Delta=cu-av\).
Then
\begin{equation}\label{eq:sczech-cqd-unified}
 \Lambda(u,v;a/c)=\mathcal F(a,c;u,\Delta)+\mathcal F(c,a;v,-\Delta).
\end{equation}
\end{theorem}

\begin{proof}
By the symmetry of the two sums, it suffices to prove
\begin{equation}\label{eq:sczech-block}
 \sum_{k=0}^{N-1}P_1\!\left(\frac{qk+\Delta}{N}\right)
       \ell\!\left(\frac{z+k}{N}\right)
 =\mathcal F(N,q;z,\Delta).
\end{equation}
Replace \(k\) by \(r=qk+h\bmod N\), \(0\le r<N\). The logarithm
becomes \(\ell(mr/N+\eta)\) and its weight is
\(r/N+\delta/N-1/2\), except that the midpoint convention adds \(1/2\)
at \(r=0\) when \(\delta=0\). The \(r/N\) part is exactly the logarithm
of the CQD. The constant part sums by the identity
\(\prod_{r=0}^{N-1}(1-\e(mr/N+\eta))=1-\e(z)\).
This gives the two cases of \eqref{eq:F-local}.
\end{proof}

\subsection{Renormalization and telescoping}

\noindent
If \(u\) or \(v\) approaches an integer, the logarithms in
\eqref{eq:sczech-w1} can diverge. Sczech subtracts their common
singular part by using the difference at \(w\) and \(1\)
\cite[\S2, (4.2)]{sczech}. For \((u,v)\in\R^2\setminus\Z^2\) put
\begin{equation}\label{eq:sczech-star}
 \Lambda_*(u,v;w)=
 \begin{cases}
 \Lambda(u,v;w),&u,v\notin\Z,\\
 \displaystyle\lim_{\substack{(s,t)\to(u,v)\\ s,t\notin\Z}}
 \bigl(\Lambda(s,t;w)-\Lambda(s,t;1)\bigr),
     &u\in\Z\ \text{or}\ v\in\Z .
 \end{cases}
\end{equation}
Sczech proves that this difference extends continuously in \((u,v)\),
and that \(\Lambda(u,v;w)\), for \(u,v\notin\Z\), is continuous in
\(w\) at each irrational \(w>0\). Here is the latter assertion for
the subtracted value.

\begin{lemma}\label{prop:sczech-ren-continuity}
If exactly one of \(u,v\) is integral, then
\(w\mapsto\Lambda_*(u,v;w)\) is continuous at every positive
irrational \(w\).
\end{lemma}

\begin{proof}
Suppose \(v\notin\Z\). After subtraction the contribution for fixed \(m\ne0\) is
\[
 \frac{\e(mu)}{m}\sum_{n\ne0}
 \bigl(\operatorname{sgn}(wm+n)-\operatorname{sgn}(m+n)\bigr)
 \frac{\e(nv)}{n}.
\]
Only integers between \(-wm\) and \(-m\), endpoints included, contribute.
On a compact subinterval of \(0<w<\infty\) these satisfy
\(|n|\asymp |m|\). The sign difference is constant in the interior.
Since \(\bigl|\sum_{n=A}^B\e(nv)\bigr|\le 2/|1-\e(v)|\), summation
by parts bounds the displayed expression, even when truncated, by
\(O_v(|m|^{-2})\), locally uniformly for \(v\notin\Z\).
Thus the sum extends to integral \(u\). At irrational \(w\), the endpoint
\(-wm\) is never an integer: for each fixed \(m\) the summand is locally
constant in \(w\). Dominated convergence now applies; the single sum
in the definition of \(\Lambda\) is absolutely convergent. The symmetry
\(\Lambda(u,v;w)=\Lambda(v,u;1/w)\) gives the other case.
\end{proof}

\noindent
Two consecutive steps share a trace: \(v_j=u_{j+1}\). If that trace
is integral, both steps require the subtraction in \eqref{eq:sczech-star}.
The two subtracted terms cancel in the limit. Thus the sum along a
finite chain still depends only on its endpoints, as follows.

\begin{lemma}\label{lem:sczech-finite-chain}
Let \(a_j,c_j\) be positive coprime integers and let \(u_j,v_j\in\R\)
for \(0\le j\le s\). Suppose that, for \(0\le j<s\) and \(b_j\in\Z\),
\[
 (a_{j+1},c_{j+1})=(c_j,b_jc_j-a_j),\qquad
 (u_{j+1},v_{j+1})=(v_j,b_jv_j-u_j).
\]
If \((u_j,v_j)\notin\Z^2\) and \(u_0,u_s\notin\Z\), then, with
\(A=c_0u_0-a_0v_0\),
\begin{equation}\label{eq:sczech-F-telescope}
 \sum_{j=0}^{s-1}\Lambda_*(u_j,v_j;a_j/c_j)
 =\mathcal F(a_0,c_0;u_0,A)-\mathcal F(a_s,c_s;u_s,A).
\end{equation}
\end{lemma}

\begin{proof}
Direct substitution gives \(c_ju_j-a_jv_j=A\). If all coordinates are nonintegral, oddness and
periodicity of \(P_1\) in \eqref{eq:sczech-block} give
\(\mathcal F(c,a;v,-A)=-\mathcal F(c,bc-a;v,A)\).
All intermediate terms in \eqref{eq:sczech-cqd-unified} cancel.

Otherwise perturb the data to
\((u_j+t a_j,v_j+t c_j)\). This keeps \(A\) fixed, hence all
Bernoulli weights fixed, and the endpoint expressions converge.
It suffices to show that the two values at \(w=1\) subtracted at
an interior integer cancel in the limit, as observed formally by
Sczech \cite[p.~141]{sczech}. Write \(p,q,r=bq-p\) for the three
consecutive perturbed coordinates, with \(q(0)\in\Z\). Then
\(p(0)\notin\Z\) and \(r(0)\equiv-p(0)\bmod\Z\). In
\(\Lambda(p,q;1)+\Lambda(q,r;1)\) the coefficient of the singular
logarithm \(\ell(q)\) is \(P_1(q-r)-P_1(p-q)=O(t)\):
the two arguments agree modulo \(\Z\) at \(t=0\), away from a jump.
Its contribution is \(O(t\log|t|)\). The other two logarithms have
the same limit and opposite limiting coefficients. Thus the reference
terms cancel and \eqref{eq:sczech-F-telescope} holds also on the axes.
\end{proof}

\subsection{From a ray class to Sczech's sum}
\label{subsec:sczech-ray-sum}

\noindent
We now construct the sum that gives \(\log X_1(\mathfrak c)\).
Each summand needs a slope \(w\) and two arguments \(u,v\).
The continued fraction supplies the slopes. To see where \(u,v\)
come from, start with an integral ideal \(\mathfrak a\) prime to
\(\frakf\) representing \(\mathfrak c\). As in
Definition~\ref{def:ray-orbit}, the lattice and the point attached
to this representative are
\[
 L=\frakf\mathfrak a^{-1},\qquad 1+L\in K/L.
\]
The lattice on which we take the continued fraction is the dual of
\(L\) for the field trace:
\[
 L^\vee=\{z\in K:\operatorname{Tr}_{K/\Q}(zL)\subseteq\Z\}
       =\frac{1}{\sqrt{D_K}}\,\mathfrak a\frakf^{-1}.
\]
Its role is to record the point \(1+L\). Indeed, for \(z\in L^\vee\)
and \(\lambda\in L\) we have
\(\operatorname{Tr}((1+\lambda)z)\equiv\operatorname{Tr}(z)\pmod{\Z}\).
Thus pairing \(1+L\) with a basis \((\alpha,\beta)\) of \(L^\vee\)
gives precisely the two arguments
\(\operatorname{Tr}(\alpha),\operatorname{Tr}(\beta)\) modulo \(\Z\).

\noindent
Choose a \(\Z\)-basis \((\alpha_0,\beta_0)\) of \(L^\vee\) with
\(\alpha_0,\beta_0\gg0\) and \(0<w_0'<1<w_0=\alpha_0/\beta_0\).
At step \(k\), record
\[
 w_k=\alpha_k/\beta_k,\qquad
 u_k=\operatorname{Tr}(\alpha_k),\qquad v_k=\operatorname{Tr}(\beta_k),
\]
and take the next basis to be
\((\alpha_{k+1},\beta_{k+1})=(\beta_k,\lceil w_k\rceil\beta_k-\alpha_k)\).
This is the minus continued-fraction step, so
\(w_{k+1}=1/(\lceil w_k\rceil-w_k)\) and \(u_{k+1}=v_k\).
In Sczech's notation \cite[(4.4)]{sczech} our order is
\((s_{k-1},s_k)=(u_k,v_k)\). The initial traces are therefore reversed
relative to the assignment preceding his Theorem~2: interchanging
the traces in \(\Lambda(u,v;w)\) also requires replacing \(w\) by \(1/w\).

\noindent
When do we stop? The slope alone need not determine the summand.
Let \(h\) be the first positive index for which
\[
 w_h=w_0,\qquad u_h\equiv u_0\pmod{\Z},\qquad
 v_h\equiv v_0\pmod{\Z}.
\]
Such an index exists: after one circuit of the reduced continued
fraction the basis has changed by a totally positive unit, which
permutes the finite orbit of \(1+L\). The integer \(h\) is the
\emph{period}; Sczech's \emph{sum} consists of the \(h\) terms
\[
 S(\mathfrak c)=\sum_{k=0}^{h-1}\Lambda_*(u_k,v_k;w_k).
\]
The two traces at a step cannot both be integers: this would imply
\(1\in(L^\vee)^\vee=L\), whereas \(L\cap\mathcal O_K=\frakf\) is
a proper ideal. Since \(u_{k+1}=v_k\), we may number the steps so
that \(u_0\notin\Z\); this cyclic permutation does not change the sum.

\noindent
To evaluate this sum by Lemma~\ref{lem:sczech-finite-chain}, we keep
its \(h\) terms and its traces fixed, and approximate only the slopes.
Put \(u=u_0,\ v=v_0\), and form the matrix for these \(h\) steps:
\[
 H=\prod_{k=0}^{h-1}\begin{pmatrix}\lceil w_k\rceil&-1\\1&0\end{pmatrix},
 \qquad H^n=\begin{pmatrix}a_n&b_n\\c_n&d_n\end{pmatrix}.
\]
The first column of \(H^n\) gives the rational starting slope for
Sczech's sum. The last column of \(H^{-n}\) gives the fraction for
the cyclic quotient:
\[
 w_{0,n}=a_n/c_n\longrightarrow w_0,\qquad
 \binom{m_n}{N_n}=H^{-n}\binom01=\binom{-b_n}{a_n},
 \qquad m_n/N_n\longrightarrow w_0'.
\]
Indeed, \(w_0\) and \(w_0'\) are the attracting fixed points of
\(H\) and \(H^{-1}\), respectively. Starting at \(w_{0,n}\), use the
same digits as before:
\(w_{k+1,n}=1/(\lceil w_k\rceil-w_{k,n})\). For each \(k<h\)
we then have \(w_{k,n}\to w_k\) as \(n\to\infty\).

\begin{theorem}\label{thm:sczech-global-corrected}
For the strict ray class \(\mathfrak c\) and its conjugate
\(\mathfrak c^\sigma\) one has
\begin{equation}\label{eq:sczech-global-corrected}
 \log X_1(\mathfrak c)=S(\mathfrak c),\qquad
 \log X_2(\mathfrak c)=S(\mathfrak c^\sigma).
\end{equation}
More precisely, with \(x=\{-v\}\), for all sufficiently large \(n\)
\begin{equation}\label{eq:sczech-rational-periodQ}
 \sum_{k=0}^{h-1}\Lambda_*(u_k,v_k;w_{k,n})
   =\log C_H(m_n,N_n;x,u)
   =S(\mathfrak c)+O(N_n^{-2}).
\end{equation}
The same formula applies to \(\mathfrak c^\sigma\) with its own data
and modulus \(\frakf^\sigma\).
\end{theorem}

\begin{proof}
We first prove the exact identity at rational slopes. The \(h\)
continued-fraction steps take the first column of \(H^n\) to that
of \(H^{n-1}\). By Lemma~\ref{lem:sczech-finite-chain} their sum
is the difference of the two \(\mathcal F\)-values at these columns.
To evaluate them write
\[
 \binom{y}{-A}=H^{-n}\binom uv
 =H^{-(n-1)}\binom{u_h}{v_h}\equiv\binom uv\pmod{\Z^2}.
\]
Thus \(y\equiv u\bmod\Z\), \(\{A\}=x\),
\(u=N_n y+m_n A\), and \(m_n c_n\equiv1\bmod N_n\).
Eventually \(0<m_n/N_n<1\), since \(m_n/N_n\to\theta=w_0'\in(0,1)\).
Hence the argument \(\eta\) in \eqref{eq:F-local} is
\(y+xm_n/N_n\). At the other endpoint the same \(y,A\) occur, with
\(n-1\) in place of \(n\) and \(u_h\equiv u\).
In the difference of the two blocks the common
\(-\ell(u)/2\) cancels, as does the midpoint term if
\(x=0\). What remains is exactly \eqref{eq:Q-nonres} for \(H\):
its input and output columns are \((m_n,N_n)\trans\) and
\((m_{n-1},N_{n-1})\trans\), and their common phase is
\(\e(xm_n+uN_n)=\e(u)\ne1\), so no zero is omitted. This proves
the first equality in \eqref{eq:sczech-rational-periodQ}.

\noindent
Now let \(n\to\infty\). There are only \(h\) terms on the left.
Each tends to its value at \(w_k\), by Sczech's continuity when
neither trace is integral and by Lemma~\ref{prop:sczech-ren-continuity}
otherwise. Thus their sum tends to \(S(\mathfrak c)\).
To identify the same limit on the cyclic side, note that \(H\) fixes
\(r=(-u,x)\trans\) modulo \(\Z^2\). One circuit of the reduced
continued fraction of the slope alone corresponds to the fundamental totally positive unit
\cite[proof of Lemma~2.1.3]{yamamoto-kronecker}. A unit power \(\gamma\) fixes both
coordinates modulo \(\Z\) precisely when
\(\operatorname{Tr}_{K/\Q}((\gamma-1)L^\vee)\subseteq\Z\), that is,
\(\gamma-1\in L\cap\mathcal O_K=\frakf\). Hence the eigenvalues of \(H\)
at \(w_0,\theta\) are \(\ep^{g(\frakf)},\ep^{-g(\frakf)}\).
Thus \(m_n/N_n=\theta+O(N_n^{-2})\) by diagonalization, and
\eqref{eq:general-rm-rate} gives
\[
 \log C_H(m_n,N_n;x,u)=\log|\mathscr S_r(H;\theta)|+O(N_n^{-2}).
\]

\noindent
It remains to recognize this cocycle value in Kopp's formula
\eqref{eq:kopp-endpoints}. The basis of \(L\) dual
to \((\alpha_0,\beta_0)\) under the trace is
\((W,-\theta W)\), where
\(W=\beta_0'/(\alpha_0\beta_0'-\alpha_0'\beta_0)\); thus
\(1=W(u-\theta v)\). To express \(1+L\) in the coordinates
of \eqref{eq:kopp-endpoints}, set
\[
 t=W(1-\theta)\gg0,\quad \xi=\theta/(1-\theta),\quad
 q=\lceil-\xi'\rceil,\quad
 M=\begin{pmatrix}1-q&q\\-1&1\end{pmatrix},\quad \omega=M\theta=\xi+q.
\]
Then \(0<\omega'<1<\omega\), \(t^{-1}L=\Z\omega+\Z\), and
\(t^{-1}=\widetilde r_2\omega-\widetilde r_1\) with
\(\widetilde r=M(-u,-v)\trans\). These identities give both the
lattice and its marked point. Since \(1-\theta>0\), the change
of basis \eqref{eq:cqd-covariance} preserves the limiting modulus.
At \(\omega\) the eigenvalue of \(MHM^{-1}\) is \(\ep^{-g(\frakf)}\),
so this is the contracting generator in \eqref{eq:kopp-endpoints}.
That formula identifies the modulus with \(X_1(\mathfrak c)\).
Finally \(X_1(\mathfrak c^\sigma)=X_2(\mathfrak c)\) by
Proposition~\ref{prop:galois-exchange}.
Consequently the sums do not depend on the chosen representative
or reduced basis, and we need not assume \(\frakf^\sigma=\frakf\).
The signature identities \eqref{eq:signature-1}--\eqref{eq:signature-2}
also give \(S(\mathfrak c\mathfrak s_1)=S(\mathfrak c)\) and
\(S((\mathfrak c\mathfrak s_2)^\sigma)=S(\mathfrak c^\sigma)\).
\end{proof}

\begingroup
\makeatletter
\let\@tocwrite\@gobbletwo
\makeatother
\section*{Declaration of AI use}
\noindent
AI tools assisted with drafting, proof review and the preparation of
SageMath checks. All mathematical ideas are due to the author, who is responsible for the final text and
its mathematical claims.

\endgroup

\end{document}